\documentclass[11pt,a4paper,reqno]{amsart}
\usepackage{amsfonts,amssymb,amsmath,amsthm,amscd,latexsym,graphicx,url, color,cite,comment,upgreek,esint,xspace}
\usepackage{a4wide}

\usepackage{xspace}
\usepackage{verbatim}
\usepackage{comment}

\usepackage[margin=1in]{geometry}

\usepackage{amsfonts,amssymb,amsmath,amsthm,amscd,latexsym,graphicx,url, color,cite,comment,upgreek,esint,xspace}
\usepackage{a4wide}

\usepackage{xspace}
\usepackage{verbatim}
\usepackage{comment}

\usepackage[margin=1in]{geometry}

\title[Boundary singularities of solutions to semilinear fractional equations]{Boundary singularities of solutions to semilinear fractional equations}

\author{Phuoc-Tai Nguyen} 
\address{Department of Mathematics and Statistics, Masaryk University, Brno, Czech Republic}
\email{ptnguyen@math.muni.cz; \, nguyenphuoctai.hcmup@gmail.com}
\author{Laurent V\'eron}
\address{Laboratoire de Math\'ematiques et Physique Th\'eorique  \\Universit\'e Fran\c{c}ois Rabelais,  Tours,  France}
\email{veronl@univ-tours.fr}

\numberwithin{equation}{section}
\newcommand{\be}{\begin{equation}}
\newcommand{\bel}[1]{\begin{equation}\label{#1}}
\newcommand{\ee}{\end{equation}}
\newcommand{\barr}{\begin{eqnarray}}
\newcommand{\earr}{\end{eqnarray}}
\newcommand{\bars}{\begin{eqnarray*}}
\newcommand{\ears}{\end{eqnarray*}}

\newtheorem{subn}{\name}

\newcommand{\bsn}[1]{\def\name{#1}\begin{subn}}
\newcommand{\esn}{\end{subn}}
\newtheorem{sub}{\name}[section]
\newcommand{\bs}{\begin{sub}}
\newcommand{\es}{\end{sub}}
\newcommand{\bsl}[1]{\begin{sub}\label{#1}}
\newcommand{\bth}[1]{\def\name{Theorem}
\begin{sub}\label{t:#1}}
\newcommand{\blemma}[1]{\def\name{Lemma}
\begin{sub}\label{l:#1}}
\newcommand{\bcor}[1]{\def\name{Corollary}
\begin{sub}\label{c:#1}}
\newcommand{\bdef}[1]{\def\name{Definition}
\begin{sub}\label{d:#1}}
\newcommand{\bprop}[1]{\def\name{Proposition}
\begin{sub}\label{p:#1}}
\newcommand{\rth}[1]{Theorem~\ref{t:#1}}
\newcommand{\rlemma}[1]{Lemma~\ref{l:#1}}
\newcommand{\rcor}[1]{Corollary~\ref{c:#1}}
\newcommand{\rdef}[1]{Definition~\ref{d:#1}}
\newcommand{\rprop}[1]{Proposition~\ref{p:#1}}
\newcommand{\BA}{\begin{array}}
\newcommand{\EA}{\end{array}}
\newcommand{\BAN}{\renewcommand{\arraystretch}{1.2}
\setlength{\arraycolsep}{2pt}\begin{array}}
\newcommand{\BAV}[2]{\renewcommand{\arraystretch}{#1}
\setlength{\arraycolsep}{#2}\begin{array}}
\newcommand{\BSA}{\begin{subarray}}
\newcommand{\ESA}{\end{subarray}}
\newcommand{\BAL}{\begin{aligned}}
\newcommand{\EAL}{\end{aligned}}
\newcommand{\BALG}{\begin{alignat}}
\newcommand{\EALG}{\end{alignat}}
\newcommand{\BALGN}{\begin{alignat*}}
\newcommand{\EALGN}{\end{alignat*}}

\newcommand{\qeda}{\hspace{10mm}\hfill $\square$}

\newcommand{\forevery}{\quad \forall}

\newcommand{\Lra}{\Longrightarrow}

\newcommand{\abs}[1]{\left |#1\right |}
\newcommand{\norm}[1]{\left \|#1\right \|}

\def\angb<#1>{\langle #1 \rangle}

\newcommand{\opname}[1]{\mbox{\rm #1}\,}
\newcommand{\supp}{\opname{supp}}
\newcommand{\dist}{\opname{dist}}
\newcommand{\myfrac}[2]{{\displaystyle \frac{#1}{#2} }}
\newcommand{\myint}[2]{{\displaystyle \int_{#1}^{#2}}}

\newcommand{\q}{\quad}
\newcommand{\qq}{\qquad}

\newcommand{\prt}{\partial}
\newcommand{\sms}{\setminus}

\newcommand{\ti}{\times}

\newcommand{\tl}{\tilde}
\newcommand{\sbs}{\subset}

\newcommand{\nind}{\noindent}

\newcommand{\ovl}{\overline}

\newcommand{\nin}{\not\in}

\def\ga{\alpha}     \def\gb{\beta}       \def\gg{\gamma}
       \def\gd{\delta}      \def\ge{\epsilon}
\def\gth{\theta}                         \def\vge{\varepsilon}
\def\gf{\phi}       \def\vgf{\varphi}    
      \def\gk{\kappa}      \def\gl{\lambda}
\def\gm{\mu}                 \def\gp{\pi}
    \def\gr{\rho}        
\def\gs{\sigma}       \def\gt{\tau}
      \def\gw{\omega}
                \def\gz{\zeta}
\def\Gg{\Gamma}     \def\Gd{\Delta}      

    \def\Gs{\Sigma}      
\def\Gw{\Omega}              

\def\CS{{\mathcal S}}      
   \def\CO{{\mathcal O}}   
\def\CA{{\mathcal A}}   \def\CB{{\mathcal B}}   
\def\CD{{\mathcal D}}   \def\CE{{\mathcal E}}   \def\CF{{\mathcal F}}
      
\def\CJ{{\mathcal J}}      \def\CL{{\mathcal L}}

\def\BBA {\mathbb A}       
   \def\BBE {\mathbb E}    
\def\BBG {\mathbb G}       
       
\def\BBM {\mathbb M}   \def\BBN {\mathbb N}    
\def\BBP {\mathbb P}   \def\BBR {\mathbb R}    
       
   \def\BBX {\mathbb X}

\def\GTM {\mathfrak M}

\def\fw{{(-\Gd)^s}}

\def\sign{\mathrm{sign\,}}
\def\tr{\mathrm{tr\,}}

\begin{document}
\maketitle \medskip

\begin{abstract}
We prove the existence of a solution of $ (-\Delta)^s u+f(u)=0$ in a smooth bounded domain $\Omega$ with  a prescribed boundary value $\mu$ in the class of  Radon measures for a large class of continuous functions $f$ satisfying a weak singularity condition expressed under an integral form. We study the existence of  a boundary trace for positive moderate solutions. In the particular case where $f(u)=u^p$ and $\mu$ is a Dirac mass, we show the existence of several critical exponents $p$. We also demonstrate the existence of several types of separable solutions of the equation $ (-\Delta)^s u+u^p=0$ in $\BBR^N_+$.

\vspace{2mm}

\end{abstract}

\tableofcontents 
\section{Introduction} \setcounter{equation}{0}
Let  $\Gw \sbs \BBR^N$ be a bounded domain with $C^2$ boundary and $s \in (0,1)$. Define the $s$-fractional Laplacian as 
$$ \fw u(x) := \lim_{\vge \to 0}(-\Gd)^s_{\vge}u(x) $$
where 
$$ (-\Gd)^s_{\vge}u(x): = a_{N,s}\int_{\BBR^N \sms B_\vge(x)}\frac{u(x)-u(y)}{|x-y|^{N+2s}}dy, \quad a_{N,s}:=\frac{\Gg(N/2+s)}{\pi^{N/2}\Gg(2-s)}s(1-s). $$ 

We denote by $G_s^\Gw$ and $M_s^\Gw$ the Green kernel and the Martin kernel of $\fw$ in $\Gw$ respectively. Denote by $\BBG_s^\Gw$ and $\BBM_s^\Gw$ the Green operator and the Martin operator (see section 2 for more details). Further, for $\gf \geq 0$, denote by $\GTM(\Gw,\gf)$ the space of Radon measures $\gt$ on $\Gw$ satisfying $\int_{\Gw}\gf d|\gt|<\infty$ and by $\GTM(\prt \Gw)$ the space of bounded Radon measures on $\prt \Gw$. Let $\rho(x)$ be the distance from $x$ to $\prt \Gw$. For $\gb>0$, set
$$ \Gw_\gb:=\{x \in \Gw: \gr(x)<\gb \},\; D_\gb:=\{x \in \Gw: \gr(x)>\gb\}, \; \Gs_\gb:=\{x \in \Gw: \gr(x)=\gb\}.$$
\bdef{trs} We say that a function $u \in L^1_{loc}(\Gw)$ possesses a $s$-boundary trace on $\prt \Gw$ if there exists a measure $\mu \in \GTM(\prt \Gw)$ such that
\bel{str} \lim_{\gb \to 0}\gb^{1-s}\int_{\Gs_\gb}|u - \BBM_s^\Gw[\gm]|dS = 0. \ee
The $s$-boundary trace of $u$ is denoted by $\tr_s(u)$. 
\es

Let $\gt \in \GTM(\Gw,\gr^s)$, $\gm \in \GTM(\prt \Gw)$ and $f \in C(\BBR)$ be a nondecreasing function with $f(0)=0$. In this paper, we study boundary singularity problem for semilinear fractional equation of the form
\bel{Np} \left\{ \BAL \fw u + f(u) &= \gt \qquad  \text{in } \Gw \\ 
\tr_s(u) &= \mu  \\
u &= 0 \qquad \text{in } \Gw^c.
\EAL \right. \ee
We denote by $\BBX_s(\Gw) \subset C(\BBR^N)$ the space of test functions $\xi$ satisfying

(i) $\supp (\xi) \sbs \bar \Gw$,

(ii) $\fw \xi(x)$ exists for all $x \in \Gw$ and $|\fw \xi(x)| \leq C$ for some $C>0$,

(iii) there exists $\vgf \in L^1(\Gw,\rho^s)$ and $\ge_0>0$ such that $|(-\Gd)_\ge^s \xi| \leq \vgf$ a.e. in $\Gw$, for all $\ge \in (0,\ge_0]$.
\bdef{defNp} Let $\gt \in \GTM(\Gw,\gr^s)$ and $\mu \in \GTM(\prt \Gw)$. A function $u$ is called a weak solution of \eqref{Np} if $u \in L^1(\Gw)$, $f(u) \in L^1(\Gw,\gr^s)$ and 
\bel{intN} \int_{\Gw} (u \fw \xi + f(u) \xi)\,dx = \int_{\Gw}\xi d\gt + \int_{\Gw} \BBM_s^\Gw[\mu] \fw \xi \,dx, \forevery \xi \in \BBX_s(\Gw). \ee
\es
The boundary value problem with measure data for semilinear elliptic equations
\bel{localpro} \left\{ \BAL -\Gd u + f(u) &= 0 \qquad  &&\text{in } \Gw \\ 
u &= \mu  &&\text{on } \prt \Gw,\\
\EAL \right. \ee
was first studied by A. Gmira and L. V\'eron in \cite{GV} and then the typical model, i.e. problem \eqref{localpro} with $f(u)=u^p$ ($p>1$), has been intensively investigated by numerous authours (see \cite{MV1,MV2,MV3,MV4,MVbook} and references therein). They proved that if $f$ is a continuous, nondecreasing function satisfying
\bel{cond-f} \int_{1}^{\infty}[f(t)-f(-t)]s^{-1-p_c}dt<\infty, \ee
where $p_c:=\frac{N+1}{N-1}$, then problem \eqref{localpro} admits a unique weak solution. In particular, when $f(u)=u^p$ with $1<p<p_c$ and $\mu=k\gd_0$ with $0 \in \prt \Gw$ and $k>0$, there exists a unique solution $u_k$ of \eqref{localpro}. It was showed \cite{MV1,MVbook} that the sequence $\{ u_k \}$ is increasing and converges to a function $u_\infty$ which is a solution of the equation in \eqref{localpro}.  
 
To our knowledge, few papers concerning boundary singularity problem for nonlinear fractional elliptic equation have been published in the literature.  The earliest works in this direction are the papers \cite{FQ,CFQ} by P. Felmer et al. which deal with the existence, nonexistence and asymptotic behavior of large solutions for equations involving fractional Laplacian. Afterwards, N. Abatangelo \cite{Ab1} presented a suitable setting for the study of fractional Laplacian equations in a measure framework and provided a fairly comprehensive description of large solutions which improve the results in \cite{FQ,CFQ}. Recently, H. Chen et al. \cite{CAHM} investigated semilinear elliptic equations involving measures concentrated on the boundary by employing approximate method. 

In the present paper, we aim to establish the existence and uniqueness of weak solutions of \eqref{Np}. To this end, we develop a theory for linear equations associated to \eqref{Np} 
\bel{Lp} \left\{ \BAL \fw u &= \gt \qquad &&\text{in } \Gw \\ 
\tr_s(u) &= \mu \qquad &&\text{} \\
u &= 0 \qquad &&\text{in } \Gw^c.
\EAL \right. \ee
Existence and uniqueness result for \eqref{Lp} is stated in the following proposition. \smallskip

\noindent \textbf{Proposition A.} \textit{Assume $s \in (\frac{1}{2},1)$. Let  $\gt \in \GTM(\Gw,\gr^s)$ and $\mu \in \GTM(\prt \Gw)$. Then problem \eqref{Lp} admits a unique weak solution. The solution is given by
\bel{RFL} u = \BBG_s^\Gw[\gt] + \BBM_s^\Gw[\mu]. \ee
Moreover, there exists a positive constant $c=c(N,s,\Gw)$ such that
\bel{esL} \norm{u}_{L^1(\Gw)} \leq c(\norm{\gt}_{\GTM(\Gw,\gr^s)}+\norm{\mu}_{\GTM(\prt \Gw)}). \ee }

This proposition allows to study semilinear equation \eqref{Np}. We first deal with the case of $L^1$ data. \smallskip

\noindent \textbf{Theorem B.} \textit{Assume $s \in (\frac{1}{2},1)$. Let $f \in C(\BBR)$ be a nondecreasing function satisfying $t f(t) \geq 0$ for every $t \in \BBR$.} 

\noindent {\sc I. Existence and uniqueness.}  \textit{For every $\gt \in L^1(\Gw,\gr^s)$ and $\mu \in L^1(\prt \Gw)$, problem \eqref{Np} admits a unique weak solution $u$. Moreover,
\bel{form} u = \BBG_s^\Gw[\gt - f(u)] + \BBM_s^\Gw[\mu] \quad \text{in } 
\Gw, \ee
\bel{estu} -\BBG_s^\Gw[\gt^-] - \BBM_s^\Gw[\gm^-] \leq u \leq  \BBG_s^\Gw[\gt^+] + \BBM_s^\Gw[\gm^+] \quad \text{in } \Gw. 
\ee}

\noindent {\sc II. Monotonicity.} \textit{The mapping $(\gt,\gm) \mapsto u$ is nondecreasing.} \medskip

\noindent \textbf{Remark.} The restriction $s \in (\frac{1}{2},1)$ is due to the fact that in this range of $s$, $\tr_s(\BBG[\gt])=0$ for every $\gt \in \GTM(\Gw,\gr^s)$ (see \rprop{trG}). We conjecture that this still holds if $s \in (0,\frac{1}{2}]$. \smallskip 

We reveal that, in measures framework, because of the interplay between the nonlocal operator $\fw$ and the nonlinearity term $f(u)$,  the analysis is much more intricate and there are 3 critical exponents
$$p_1^*:=\frac{N+2s}{N}, \quad p_2^*:=\frac{N+s}{N-s}, \quad p_3^*:=\frac{N}{N-2s}.$$
This yields substantial new difficulties and leads to disclose new types of results. The new aspects are both on the technical side and on the one of the new phenomena observed. \medskip

\noindent \textbf{Theorem C.} \textit{Assume $s \in (\frac{1}{2},1)$. Let $f \in C(\BBR)$ be a nondecreasing function, $t f(t) \geq 0$ for every $t \in \BBR$ and}
\bel{subcri} \int_1^\infty [f(s)-f(-s)]s^{-1-p_2^*}ds < \infty. \ee
\noindent {\sc I. Existence and Uniqueness.} \textit{For every $\gt \in \GTM(\Gw,\gr^s)$ and $\mu \in \GTM(\prt \Gw)$ there exists a unique weak solution of \eqref{Np}. This solution satisfies \eqref{form} and \eqref{estu}. Moreover, the mapping $(\gt,\gm) \mapsto u$ is nondecreasing.} \smallskip

\noindent {\sc II. Stability.} \textit{Assume $\{\gt_n\} \sbs \GTM(\Gw,\gr^s)$ converges weakly to $\gt \in \GTM(\Gw,\gr^s)$ and $\{\gm_n\} \sbs \GTM(\prt \Gw)$ converges weakly to $\gm \in \GTM(\prt \Gw)$. Let $u$ and $u_n$ be the unique weak solutions of \eqref{Np} with data $(\gt,\gm)$ and $(\gt_n,\gm_n)$ respectively. Then $u_n \to u$ in $L^1(\Gw)$ and $f(u_n) \to f(u)$ in $L^p(\Gw,\gr^s)$.
} \smallskip

If $\mu$ is a Dirac mass concentrated at a point on $\prt \Gw$, we obtain the behavior of the solution near that boundary point. \smallskip
 
\noindent \textbf{Theorem D.} \textit{Under the assumption of Theorem C, let $z \in \prt \Gw$, $k>0$ and $u_{z,k}^\Gw$ be the unique weak solution of 
\bel{Npd} \left\{ \BAL \fw u + f(u) &= 0 \qquad \text{in } \Gw \\ 
\tr_s(u) &= k\gd_z  \\
u &= 0 \qquad \text{in } \Gw^c.
\EAL \right. \ee
Then 
\bel{ukz} \lim_{\Gw \ni x \to z}\frac{u_{z,k}^\Gw(x)}{M_s^\Gw(x,z)}= k. \ee
}

We next assume that $0 \in \prt \Gw$. Let $0<p<p_2^*$ and denote by $u_k^\Gw$ the unique weak solution of 
\bel{Nppd} \left\{ \BA{rll} \fw u + u^p &= 0 \qquad &\text{in } \Gw \\ 
\tr_s(u) &= k\gd_0  \\
u &= 0 \qquad &\text{in } \Gw^c.
\EA \right. \ee

By Theorem C, $u_k^\Gw \leq kM_s^\Gw(\cdot,0)$ and $k \mapsto u_k^\Gw$ is increasing. Therefore, it is natural to investigate $\lim_{k \to \infty}u_k^\Gw$.  This is accomplishable thanks to the study of separable solutions of 
\bel{I1} \left\{ \BA{rll} \fw u + u^p &= 0 \qquad &\text{in } \BBR^N_+ \\ 
u &= 0 \qquad &\text{in } \overline{\BBR^N_-}
\EA \right. \ee
with $p>1$. Denote by 
$$S^{N-1}:=\left\{\gs=(\cos\gf\,\gs',\sin\gf):\gs'\in S^{N-2}, -\tfrac\gp2\leq\gf\leq \tfrac\gp2\right\} $$
the unit sphere in $\BBR^N$ and by $S_+^{N-1}:=S^{N-1} \cap \BBR_+^N$ the upper hemisphere. 
Writing separable solution under the form $u(x)=u(r,\gs)=r^{-\frac{2s}{p-1}}\gw(\gs)$, with $r>0$ and 
$\gs\in S^{N-1}_+$, we obtain that $\gw$ satisfies 
\bel{I2}\left\{\BA {lll}
\CA_s\gw-\CL_{s,\frac{2s}{p-1}}\gw+\gw^p=0\qquad&\text{in }\;S_+^{N-1}\\[2mm]
\phantom{\CA_s\gw-\CL_{s,\frac{2s}{p-1}}+\gw^p}
\gw=0&\text{in }\;\overline{S_-^{N-1}},
\EA\right.\ee
where $\CA_s$ is a nonlocal operator naturally associated to the $s$-fractional Laplace-Beltrami operator and 
$\CL_{s,\frac{2s}{p-1}}$ is a linear integral operator with kernel. In analyzing the spectral properties of  $\CA_s$ we prove \medskip

\noindent \textbf{Theorem E.} \textit{Let $N\geq 2$, $s\in (0,1)$ and $p>p_1^*$.} \smallskip

\nind \textit{I- If $p_2^* \leq p<p_3^*$ there exists no positive solution of \eqref{I2} belonging to $W^{s,2}_0(S^{N-1}_+)$.}
\smallskip

\nind \textit{II- If $p_1^*<p<p_2^*$ there exists a unique positive solution  $\gw^*\in W^{s,2}_0(S^{N-1}_+)$ of \eqref{I2}.} \medskip

As a consequence of this result we obtain the behavior of $u_k^\Gw$ when $k\to\infty$.\smallskip

\noindent \textbf{Theorem F} \textit{Assume $s \in (\frac{1}{2},1)$. Let $\Gw=\BBR^N_+$ or $\Gw$ be a bounded domain with $C^2$ boundary containing $0$.} 

\noindent \textit{I- If $p \in (p_1^*,p_2^*)$ then $u_\infty^\Gw:=\lim_{k \to 0}u_k^\Gw$ is a positive solution of 
\bel{eqa} \left\{ \BA{rll} \fw u + u^p &= 0 \qquad &\text{in } \Gw \\ 
u &= 0 \qquad &\text{in } \Gw^c.
\EA \right. \ee}

\textit{(i) If $\Gw=\BBR_+^N$ then 
$$ u_\infty^{\BBR_+^N}(x)=|x|^{-\frac{2s}{p-1}}\gw^*(\gs), \quad \text{with } \gs=\frac{x}{|x|} \qquad \forall x \in \BBR_+^N.$$}

\textit{(ii) If $\Gw$ is a bounded $C^2$ domain with $\prt \Gw$ containing $0$ then 
\bel{min1} 
\lim_{\tiny \BA{c}\Gw \ni x\to 0\\
\frac{x}{|x|}=\gs\in S^{N-1}_+
\EA}|x|^{\frac{2s}{p-1}}u_{\infty}^{\Gw}(x)=\gw^*(\gs),
\ee
 locally uniformly on $S^{N-1}_+$. In particular, there exists a positive constant $c$ depending on $N$, $s$, $p$ and the $C^2$ norm of $\prt\Gw$ such that
\bel{uinf2} c^{-1}\gr(x)^s|x|^{-\frac{(p+1)s}{p-1}} \leq u_\infty^\Gw(x) \leq c\gr(x)^s|x|^{-\frac{(p+1)s}{p-1}} \forevery x \in \Gw.
\ee
II- Assume $p \in (0,p_1^*]$. Then $\lim_{k \to \infty}u_k^\Gw=\infty$ in $\Gw$.
} \medskip

The main ingredients of the present study: estimates on Green kernel and Martin kernel, theory for linear fractional equations in connection with the notion $s-$boundary trace as mentioned above, similarity transformation and the study of equation \eqref{I2}. 

The paper is organized as follows. In Section 2, we present important properties of $s$-boundary trace and prove Proposition A. Theorems B,C,D and F are obtained in Section 3. Finally, in Appendix, we discuss separable solutions of \eqref{I1} and demonstrate Theorem E.
 
\section{Linear problems}
Throughout the present paper, we denote by $c,c',c_1,c_2,C,...$ positive constants that may vary from line to line. If necessary, the dependence of these constants will be made precise.
\subsection{$s$-harmonic functions}
We first recall the definition of \emph{$s$-harmonic functions} (see \cite[page 46]{Bog1}, \cite[page 230]{Bog2}, \cite[page 20]{BBK}). Denote by $(X_t,P^x)$ the standard rotation invariant $2s$-stable L\'evy process in $\BBR^N$ (i.e. homogeneous with independent increments) with characteristic function 
$$E^0 e^{i\xi X_t}=e^{-t|\xi|^{2s}}, \quad \xi \in \BBR^N, t\geq 0. $$
Denote by $E^x$ the expectation with respect to the distribution $P^x$ of the process starting from $x \in \BBR^N$. We assume that sample paths of $X_t$ are right-continuous and have left-hand limits a.s. The process $(X_t)$ is Markov with transition probabilities given by 
$$ P_t(x,A)=P^x(X_t \in A)= \mu_t(A-x),$$
where $\mu_t$ is the one-dimensional distribution of $X_t$ with respect to $P^0$. It is well known that $\fw$ is  the generator of the process $(X_t,P^x)$.

For each Borel set $D \sbs \BBR^N$, set $t_D:=\inf\{ t \geq 0: X_t \nin D \}$, i.e. $t_D$ is the first exit time from $D$. If $D$ is bounded then $t_D<\infty$ a.s.  Denote
$$ E^x u(X_{t_D})=E^x\{ u(X_{t_D}): t_D<\infty \}. $$
\bdef{har1} Let $u$ be a Borel measurable function in $\BBR^N$. We say that $u$ is $s$-harmonic in $\Gw$ if for every bounded open set $D \Subset \Gw$,  
$$ u(x)=E^xu(X_{t_D}),  \quad x \in D.$$ 
We say that $u$ is singular $s$-harmonic in $\Gw$ if $u$ is $s$-harmonic and $u=0$ in $\Gw^c$.
\es
Put 
$$ \CD_s:=\left\{ u: \BBR^N \mapsto \BBR: \text{Borel measurable such that} \int_{\BBR^N}\frac{|u(x)|}{(1+ |x|)^{N+2s}} \right\}. $$
The following result follows from \cite[Corollary 3.10 and Theorem 3.12]{BB} and \cite[page 20]{BBK} (see also \cite{Kw}).
\bprop{harmonic} Let $u \in \CD_s$. 

\emph{(i)} $u$ is $s$-harmonic in $\Gw$ if and only if $\fw u = 0$ in $\Gw$ in the sense of distributions. 

\emph{(ii)} $u$ is singular $s$-harmonic  in $\Gw$ if and only if $u$ is  $s$-harmonic in $\Gw$ and $u=0$ in $\Gw^c$. 
\es

\subsection{Green kernel, Poisson kernel and Martin kernel}

In what follows the notation $f \sim g$ means: there exists a positive constant $c$ such that $c^{-1}f < g < cf$ in the domain of the two functions or in a specified  subset of this domain.

Denote by $G_s^\Gw$ the Green kernel of $\fw$ in $\Gw$. Namely, for every $y \in \Gw$, 
$$ \left\{ \BAL \fw G_s^\Gw(\cdot,y) &=\gd_y \quad &&\text{in } \Gw \\
G_s^\Gw(\cdot,y) &= 0 \quad &&\text{in } \Gw^c, \EAL \right. $$
where $\gd_y$ is the Dirac mass at $y$. By combining \cite[Lemma 3.2]{Ab1} and \cite[Corollary 1.3]{CS1}), we get
\bprop{PropG}
(i) $G_s^\Gw$ is in continuous, symmetric, positive in $\{ (x,y) \in \Gw \times \Gw: x \neq y\}$ and $G_s^\Gw(x,y)=0$ if $x$ or $y$ belongs to $\Gw^c$.  

(ii) $\fw G_s^\Gw(x,\cdot) \in L^1(\Gw^c)$ for every $x \in \Gw$ and $\fw G_s^\Gw(x,y) \leq 0$ for every $x \in \Gw$ and $y \in \Gw^c$.

(iii) There holds
\bel{Ga} G_s^\Gw(x,y) \sim
\min\left\{\abs{x-y}^{2s-N},\gr(x)^{s}\gr(y)^{s}\abs{x-y}^{-N}\right\} \quad \forall (x,y) \in \Gw \times \Gw, x \neq y. \ee
The similarity constant in the above estimate depends only on $\Gw$ and $s$. 
\es

Denote by $\BBG_s^\Gw$ the associated Green operator 
$$ \BBG_s^\Gw[\gt] = \int_{\Gw}G_s^\Gw(\cdot,y)d\gt(y) \qquad \tau \in \GTM(\Gw,\rho^s). $$
Put
\bel{k} k_{s,\gg}:=\left\{ \BA{lll} &p_3^* \qquad &\text{if } \gg \in [0,\frac{N-2s}{N}s) \\[3mm]
&\frac{N+s}{N-2s+\gg} \qquad &\text{if } \gg \in [\frac{N-2s}{N}s,s]. \EA \right. \ee

H. Chen and L. V\'eron obtained the following estimate for Green operator \cite[Proposition 2.3 and Proposition 2.6]{CV1}. 
\blemma{G} Assume $\gg \in [0,s]$ and $k_{s,\gg}$ be  as in \eqref{k}. 

\emph{(i)} There exists a constant $c=c(N,s,\gg,\Gw)>0$ such that
\bel{estG1} \norm{\BBG^\Gw_s[\gt]}_{M^{k_{s,\gg}}(\Gw,\gr^s)}
\leq c\norm{\gt}_{\GTM(\Gw,\gr^\gg)}\q \forall \gt \in \GTM(\Gw,\gr^\gg).\ee

\emph{(ii)} Assume $\{ \gt_n \} \sbs \GTM(\Gw,\gr^\gg)$ converges weakly to $\gt \in \GTM(\Gw,\gr^\gg)$. Then $\BBG_s^\Gw[\gt_n] \to \BBG_s^\Gw[\gt]$ in $L^p(\Gw,\gr^s)$ for any $p \in [1,k_{s,\gg})$. 
\es

Let $P_s^\Gw$ be the Poisson kernel of $\fw$  defined by (see \cite{BKK})
$$ P_s^\Gw(x,y):=-a_{N,-s}\int_\Gw \frac{G_s^\Gw(x,z)}{|z-y|^{N+2s}}dz, \forevery x \in \Gw, y \in {\ovl \Gw}^c. $$
The relation between $P_s^\Gw$ and $G_s^\Gw$ is expressed in \cite[Proposition 2]{Ab1} (see also \cite[Theorem 1.4]{CS1}, \cite[Lemma 2]{Bog2}, \cite[Theorem 1.5]{CS1}).
\bprop{propP}
 (i) $P_s^\Gw(x,y)=-\fw G_s^\Gw(x,y)$ for every $x \in \Gw$ and $y \in {\ovl \Gw}^c$. Moreover, $P_s^\Gw$ is continuous in $\Gw \times {\ovl \Gw}^c$.

(ii) There holds 
\bel{P} P_s^\Gw(x,y) \sim \frac{\gr(x)^s}{\gr(y)^s(1+ \gr(y))^s}\frac{1}{|x-y|^N}, \forevery x \in \Gw, y \in {\ovl \Gw}^c. \ee
The similarity constant in the above estimate depends only on $\Gw$ and $s$.
\es Denote by $\BBP_s^\Gw$ the corresponding operator defined by
$$ \BBP_s^\Gw[\nu](x)=\int_{{\ovl \Gw}^c}P_s^\Gw(x,y)d\nu(y), \quad \nu \in \GTM({\ovl \Gw}^c). $$
 
Fix a reference point $x_0 \in \Gw$ and denote by $M_s^\Gw$ the Martin kernel of $\fw$ in $\Gw$, i.e. 
$$ M_s^\Gw(x,z)=\lim_{\Gw \ni y \to z}\frac{G_s^\Gw(x,y)}{G_s^\Gw(x_0,y)}, \forevery x \in \BBR^N, z \in \prt \Gw. $$
By \cite[Theorem 3.6]{CS2}, the Martin boundary of $\Gw$ can be identified with the Euclidean boundary $\prt \Gw$. Denote by $\BBM_s^\Gw$ the associated Martin operator 
$$ \BBM_s^\Gw[\mu] = \int_{\prt \Gw}M_s^\Gw(\cdot,z)d\mu(z), \quad \mu \in \GTM(\prt \Gw). $$
The next result \cite{Bog2,CS2} is important in the study of $s$-harmonic functions, which give a unique presentation of $s-$harmonic functions in terms of Martin kernel.   
\bprop{martin}
\emph{(i)} The mapping $(x,z) \mapsto M_s^\Gw(x,z)$ is continuous on $\Gw \times \prt \Gw$. For any $z \in \prt \Gw$, the function $M_s^\Gw(.,z)$ is singular $s$-harmonic in $\Gw$ with  $M_s^\Gw(x_0,z)=1$. Moreover, if $z,z' \in \prt \Gw$, $z \neq z'$ then $\lim_{x \to z'}M_s^\Gw(x,z)=0$. 

\emph{(ii)} There holds 
\bel{estMar} M_s^\Gw(x,z) \sim \rho(x)^s|x-z|^{-N} \quad \forall x \in \Gw, z \in \prt \Gw.  \ee
The similarity constant in the above estimate depends only on $\Gw$ and $s$.

\emph{(iii)} For every $\mu \in \GTM^+(\prt \Gw)$ the function $\BBM_s^\Gw[\mu]$ is singular $s$-harmonic in $\Gw$ with $u(x_0)=\mu(\BBR^N)$. Conversely, if $u$ is a nonnegative singular $s$-harmonic function in $\Gw$ then there exists a unique $\mu \in \GTM^+(\prt \Gw)$ such that $u=\BBM_s^\Gw[\mu]$ in $\BBR^N$. 

\emph{(iv)} If $u$ is a nonnegative $s$-harmonic function in $\Gw$ then there exists a unique $\mu \in \GTM^+(\prt \Gw)$ such that 
$$u(x) = \BBM_s^\Gw[\mu](x) + \BBP_s^\Gw[u](x) \forevery x \in \Gw.$$ 
\es  

\blemma{M} \emph{(i)} There exists a constant $c=c(N,s,\gg, \Gw)$ such that
\bel{estM}
\norm{\BBM_s^\Gw[\gm]}_{M^{\frac{N+\gg}{N-s}}(\Gw,\gr^\gg)} \leq
c\norm{\gm}_{\GTM(\prt \Gw)},\q \forall \gm \in \GTM(\prt \Gw),\;  \gg>-s. \ee

\emph{(ii)} If $\{ \gm_n \} \sbs \GTM(\prt \Gw)$ converges weakly to $\mu \in \GTM(\prt \Gw)$ then $\BBM_s^\Gw[\gm_n] \to \BBM_s^\Gw[\gm]$ in $L^p(\Gw,\gr^\gg)$ for every $1 \leq p<\frac{N+\gg}{N-s}$.  
\es
\proof (i) By using \eqref{estMar} and a similar argument as in the proof of \cite[Theorem 2.5]{BVi}, we obtain \eqref{estM}. 
\smallskip

(ii) By combining the fact that $M_s^\Gw(x,z)=0$ for every $x \in \Gw^c$, $z \in \prt \Gw$ and \rprop{martin} (i) we deduce that for every $x \in \BBR^N$, $M_s^\Gw(x,\cdot) \in C(\prt \Gw)$. It follows that $\BBM_s^\Gw[\gm_n] \to \BBM_s^\Gw[\gm]$ everywhere in $\Gw$. Due to (i) and the Holder inequality, we deduce that, for any $1 \leq p \leq \frac{N+\gg}{N-s}$, $\{\BBM_s^\Gw[\gm_n]\}$ is uniformly integrable with respect to $\gr^\gg dx$. By invoking Vitali's theorem, we obtain the convergence in $L^p(\Gw,\gr^\gg)$. \qeda

\subsection{Boundary trace}
We recall that, for $\gb>0$, 
$$ \Gw_\gb:=\{x \in \Gw: \gr(x)<\gb \},\; D_\gb:=\{x \in \Gw: \gr(x)>\gb\}, \; \Gs_\gb:=\{x \in \Gw: \gr(x)=\gb\}.$$
The following geometric property of $C^2$ domains can be found in \cite{MVbook}.
\bprop{beta0} There exists $\gb_0>0$ such that

(i) For every point $x \in \ovl \Gw_{\gb_0}$, there exists a unique
point $z_x \in \prt \Gw$ such that $|x -z_x|=\gr(x)$. This implies
$x=z_x - \gr(x){\bf n}_{z_x}$.

(ii) The mappings $x \mapsto \gr(x)$ and $x \mapsto z_x$ belong to
$C^2(\ovl \Gw_{\gb_0})$ and $C^1(\ovl \Gw_{\gb_0})$ respectively. Furthermore,
$\lim_{x \to z_x}\nabla \gr(x) = - {\bf n}_{z_x}$.

\es
\bprop{extrace} Assume $s \in (0,1)$. Then there exist positive constants $c=c(N,\Gw,s)$ such that, for every
$\gb\in (0,\gb_0)$,
\bel{L1}
c^{-1}\leq \gb^{1-s}\int_{\Gs_\gb}M_s^\Gw(x,y)dS(x) \leq c
\q\forall y\in \prt\Gw.
\ee
\es
\proof For $r_0>0$ fixed, by \eqref{estMar},
\bel{L3}
\int_{\Gs_\gb \sms B_{r_0}(y)}M_s^\Gw(x,y)dS(x)\leq
c_1\gb^{s},
\ee
which implies
\bel{L2}
\lim_{\gb\to 0} \int_{\Gs_\gb\sms
	B_{r_0}(y)}M_s^\Gw(x,y)dS(x)=0 \q\forall y\in \prt\Gw.
\ee
Note that for $r_0$ fixed, the rate of convergence is independent of $y$.

In order to prove \eqref{L1} we may assume that the coordinates are
placed so that $y=0$ and the tangent hyperplane to $\prt\Gw$ at $0$ is $x_N=0$ with the $x_N$ axis pointing into the domain. For $x\in \BBR^N$ put 
$x'=(x_1,\cdots,x_{N-1})$. Pick $r_0\in (0,\gb_0)$  sufficiently small
(depending only on the $C^2$ characteristic of $\Gw$) so that
$$\frac{1}{2}(|x'|^2+\gr(x)^2)\leq |x|^2\q \forall x\in \Gw\cap B_{r_0}(0).$$
Hence if $x\in \Gs_{\gb}\cap B_{r_0}(0)$ then $\frac{1}{4}( |x'| +\gb)\leq |x|$. Combining this inequality and \eqref{estMar} leads to
\begin{equation*}\BAL
\int_{\Gs_{\gb}\cap B_{r_0}(0)} M_s^\Gw(x,0)dS(x)&\leq c_2\gb^{s}\int_{\Gs_{\gb,0}} (|x'|
+\gb)^{-N}dS(x)\\
&\leq c_2\gb^{s}\int_{|x'|<r_0} (|x'| +\gb)^{-N}dx'\\
&=c_3\gb^{s-1}.
\EAL\end{equation*}
Therefore, for $\gb<r_0$,
\bel{L4}
\gb^{1-s}\int_{\Gs_{\gb}\cap B_{r_0}(0)} M_s^\Gw(x,0)dS(x) \leq c_4.
\ee
By combining estimates \eqref{L3} and \eqref{L4}, we obtain the second estimate in \eqref{L1}. The first estimate in \eqref{L1} follows from \eqref{estMar}. \qeda\smallskip

As a consequence, we get the following estimates.
\bcor{2sides} Assume $s \in (0,1)$. For every $\gm \in \GTM^+(\prt \Gw)$ and $\gb \in (0,\gb_0)$, there holds
\bel{2side}
\BAL c^{-1} \norm{\gm}_{\GTM(\prt \Gw)}  
\leq \gb^{1-s}\int_{\Gs_\gb}\BBM_s^\Gw[\gm]dS 
\leq c\norm{\gm}_{\GTM(\prt \Gw)}, \EAL
\ee
with $c$ is as in \eqref{L1}.
\es

\bprop{trG} Assume $s \in (\frac{1}{2},1)$. 
Then there exists a constant $c=c(s,N,\Gw)$ such that for any $\tau\in \GTM(\Gw,\gr^s)$ and any $0<\gb<\gb_0$,
\begin{equation}\label{G1}
\gb^{1-s} \int_{\Gs_\gb}\BBG_s^\Gw[\tau]dS\leq c\int_\Gw\gr^{s}d|\tau|.
\end{equation}
Moreover,
\begin{equation}\label{G2}
\lim_{\gb\to 0}\gb^{1-s}\int_{\Gs_\gb}\BBG_s^\Gw[\tau]dS=0.
\end{equation}
\es
\proof Without loss of generality, we may assume that $\tau>0$. Denote $v:=\BBG_s^\Gw[\tau]$. We first prove \eqref{G1}. By Fubini's theorem and \eqref{estMar},
\[\BAL
\int_{\Gs_\gb}v(x)dS(x)  &\leq c_5\Big(\int_\Gw\int_{\Gs_\gb\cap B_{\frac{\gb}{2}}(y)}|x-y|^{2s -N}
dS(x)\,d\tau(y)\\
&+\gb^{s}\int_\Gw\int_{\Gs_\gb\sms B_{\frac{\gb}{2}}(y)}|x-y|^{-N}
dS(x)\,\gr(y)^{s}d\tau(y)\Big)\\
&:=I_{1,\gb}+I_{2,\gb}.
\EAL\]
Note that, if $x\in\Gs_\gb \cap B_{\frac{\gb}{2}}(y)$ then $\gb/2\leq \gr(y)\leq 3\gb/2$. Therefore
$$\BAL
\gb^{1-s}I_{1,\gb} &\leq c_{6} \gb^{1-2s} \int_{\Gs_\gb\cap B_{\frac{\gb}{2}}(y)}|x-y|^{2s -N}dS(x)\int_\Gw \gr(y)^s\,
d\tau(y)\\
&\leq c_{6}\gb^{1-2s}\int_0^{\gb/2}r^{2s -N}r^{N-2}dr\,\int_\Gw \gr(y)^s\, d\tau(y)\\
&\leq c_{7}\int_\Gw \gr(y)^s\, d\tau(y),
\EAL$$
where the last inequality holds since $s>\frac{1}{2}$. On the other hands, we have
$$I_{2,\gb}\leq c_{7}\gb^s\int_{\gb/2}^\infty r^{-N}r^{N-2}dr\int_\Gw \gr(y)^s\,d\tau(y)=c_{8}\gb^{s-1}\int_\Gw \gr(y)^s\,d\tau(y).$$
Combining the above estimates, we obtain \eqref{G1}.

Next we demonstrate \eqref{G2}. Given $\ge\in (0,\norm{\tau}_{\GTM(\Gw,\gr^s)})$ and $\gb_1\in (0,\gb_0)$ put $\tau_1=\tau\chi_{_{\bar
D_{\gb_1}}}$ and $\tau_2=\tau \chi_{_{
\Gw_{\gb_1}} } $. We can choose $\gb_1=\gb_1(\ge)$ such that
\begin{equation}\label{tau2}
 \int_{\Gw_{\gb_1}}\gr(y)^s\,d\tau(y)\leq\ge.
\end{equation}
 Thus the choice of $\gb_1$ depends on the rate at which $\int_{\Gw_\gb}\gr^s\,d\tau $ tends to zero as  $\gb\to 0$.

Put $v_i:=\BBG_s^\Gw[\tau_i]$.
Then, for $0<\gb<\gb_1/2$,
$$\int_{\Gs_\gb}v_1(x)\,dS(x)\leq c_{9}\gb^{s}\gb_1^{-N}\int_\Gw \gr(y)^{s}d\tau_1(y),$$
which yields
\bel{Gv1}
   \lim_{\gb\to 0}\gb^{1-s} \int_{\Gs_\gb}v_1(x)\,dS(x)=0.
\ee
On the other hand, due to \eqref{G1},
\bel{Gv2}
   \gb^{1-s} \int_{\Gs_\gb}v_2\,dS \leq c_{10} \int_{\Gw}\rho^s d \gt_2 \leq c_{11}\ge \forevery \gb<\gb_0.
\end{equation}
From \eqref{Gv1} and \eqref{Gv2}, we obtain \eqref{G2}. \qeda \medskip

\blemma{sub} Assume $s \in (\frac{1}{2},1)$. Let $u, w \in \CD_s$ be two nonnegative functions satisfying
\bel{sub} \left\{ \BAL   \fw u &\leq \,  0 \, \leq \fw w \qquad &&\text{in } \Gw, \\
u &=\,0 \qquad &&\text{in } \Gw^c.  \EAL \right. \ee
If $u \leq w$ in $\BBR^N$ then $\fw u \in \GTM(\Gw,\gr^s)$ and there exists a measure $\mu \in \GTM^+(\prt \Gw)$ such that
\bel{qw} \lim_{\gb \to 0}\gb^{1-s}\int_{\Gs_\gb}|u - \BBM_s^\Gw[\gm]|dS = 0. \ee
Moreover, if $\mu=0$ then $u=0$.
\es
\proof By the assumption, there exists a nonnegative Radon measure $\tau$ on $\Gw$ such that $\fw u = - \tau$. 

We first prove that $\gt \in \GTM^+(\Gw,\gr^s)$. Define 
\bel{tlM} \tl M_s^\Gw(x,z):=\lim_{\Gw \ni y \to z}\frac{G_s^\Gw(x,y)}{\gr(y)^s}.
\ee
By \cite[page 5547]{Ab1}, there holds is a positive constant $c=c(\Gw,s)$ such that
\bel{esttlM}  \tl M_s^\Gw(x,z) \sim \rho(x)^s|x-z|^{-N}, \forevery x \in \Gw, z \in \prt \Gw, \ee
where the similarity constant depends only on $\Gw$ and $s$. This follows
\bel{esttlM} \BAL c_{12}^{-1} < c_{13}^{-1}\int_{\prt \Gw}\rho(x)|x-z|^{-N}dS(z)& \\
\leq  \gr(x)^{1-s}&\int_{\prt \Gw}\tl M_s^\Gw(x,z)dS(z) \\
&\leq c_{13}\int_{\prt \Gw}\rho(x)|x-z|^{-N}dS(z) <c_{12} \forevery x \in \Gw. \EAL \ee 
We define
$$ \BBE_s^\Gw[u](z):=\lim_{\Gw \ni x \to z} \frac{u(x)}{\int_{\prt \Gw}\tl M_s^\Gw(x,y)dS(y)} \quad z \in \prt \Gw. 
$$
For any $\gb \in (0,\gb_0)$, denote by $\gt_\gb$ the restriction of $\gt$ to $D_\gb$ and by $v_\gb$ the restriction of $u$ on $\Gs_\gb$. By \cite[Theorem 1.4]{Ab1}, there exists a unique solution $v_\gb$ of
$$ \left\{ \BAL \fw v_\gb &= -\gt_\gb \qquad &&\text{in } D_\gb\\ 
\BBE_s^{D_\gb}[v_\gb] &= 0 \qquad &&\text{on } \Gs_\gb \\
v_\gb &= u|_{D_\gb^c} \qquad &&\text{in }  D_\gb^c.
\EAL \right. $$
Moreover, the solution can be written as
\bel{vb} v_\gb+\BBG_s^{D_\gb}[\gt_\gb] = \BBP_s^{D_\gb}[u|_{ D_\gb^c}] \quad \text{in } D_\gb. \ee
By the maximum principle \cite[Lemma 3.9]{Ab1}, $v_\gb=u$ and $\BBP_s^{D_\gb}[u|_{D_\gb^c}] \leq w$ a.e. in $\BBR^N$. This, together with \eqref{vb}, implies that $\BBG_s^{D_\gb}[\gt_\gb] \leq w$ in $D_\gb$. Letting $\gb \to 0$ yields $\BBG_s^\Gw[\gt]<\infty$. For fixed $x_0 \in \Gw$, by \eqref{Ga}, $G_s^\Gw(x_0,y) > c\gr(y)^s$ for every $y \in \Gw$. Hence the finiteness of $\BBG_s^\Gw[\gt]$ implies that $\gt \in \GTM^+(\Gw,\gr^s)$. 

We next show that there exists a measure $\mu \in \GTM^+(\prt \Gw)$ such that \eqref{qw} holds. Put $v=u + \BBG_s^\Gw[\tau]$ then $v$ is a nonnegative singular $s$-harmonic in $\Gw$ due to the fact that $\BBG_s^\Gw[\tau]=0$ in $\Gw^c$. By \rprop{harmonic} and \rprop{martin} (iii), there exists $\mu \in \GTM^+(\prt \Gw)$ such that $v=\BBM_s^\Gw[\mu]$ in $\BBR^N$. By \rprop{trG}, we obtain \eqref{qw}. If $\gm=0$ then $v=0$ and thus $u=0$. \qeda

\bdef{trs} A function $u$ possesses a $s$-boundary trace on $\prt \Gw$ if there exists a measure $\mu \in \GTM(\prt \Gw)$ such that
\bel{str} \lim_{\gb \to 0}\gb^{1-s}\int_{\Gs_\gb}|u - \BBM_s^\Gw[\gm]|dS = 0. \ee
The $s$-boundary trace of $u$ is denoted noted by $\tr_s(u)$. 
\es

\noindent \textbf{\bf Remark.} (i) The notation of $s$-boundary trace is well defined. Indeed, suppose that $\mu$ and $\mu'$ satisfy \eqref{str}. Put $v=(\BBM_s^\Gw[\mu-\mu'])^+$. Clearly $v \leq \BBM_s^\Gw[|\mu|+|\mu'|]$, $v=0$ in $\Gw^c$ and $\lim_{\gb \to 0}\gb^{1-s}\int_{\Gs_\gb}|v|dS= 0$. By Kato's inequality \cite[Theorem 1.2]{CaSi}, $\fw v \leq 0$ in $\Gw$. Therefore, we deduce $v \equiv 0$ from \rlemma{sub}. This implies $\BBM_s^\Gw[\mu-\mu'] \leq 0$. By permuting the role of $\mu$ and $\mu'$, we obtain  $\BBM_s^\Gw[\mu-\mu'] \geq 0$. Thus $\mu=\mu'$.

(ii) It is clear that for every $\mu \in \GTM(\prt \Gw)$, $\tr_s(\BBM_s^\Gw[\mu])=\mu$. Moreover, if $s>\frac{1}{2}$, by \rprop{trG}, for every $\gt \in \GTM(\Gw,\gr^s)$, $\tr_s(\BBG_s^\Gw[\gt])=0$.

(iii) This kind of boundary trace was first introduced by P.-T. Nguyen and M. Marcus  \cite{MN} in order to investigate semilinear elliptic equations with Hardy potential. In the present paper we prove that it is still an effective tools in the study of nonlocal fractional elliptic equations.
\subsection{Weak solutions of linear problems}

\bdef{defLp} Let $\gt \in \GTM(\Gw,\gr^s)$ and $\mu \in \GTM(\prt \Gw)$. A function $u$ is called a weak solution of \eqref{Lp} if $u \in L^1(\Gw)$ and 
\bel{intL} \int_{\Gw} u \fw \xi \,dx = \int_{\Gw}\xi d\gt + \int_{\Gw} \BBM_s^\Gw[\mu] \fw \xi \,dx, \forevery \xi \in \BBX_s(\Gw). \ee
\es
\noindent \textbf{Proof of Proposition A.} The uniqueness follows from \cite[Proposition 2.4]{CV1}. Let $u$ be as in \eqref{RFL}. By \cite{CV1},
$$ \int_{\Gw} (u-\BBM_s^\Gw[\mu]) \fw \xi \, dx = \int_{\Gw} \BBG_s^\Gw[\gt] \fw \xi \, dx = \int_{\Gw}\xi d\gt \quad \forall \xi \in \BBX_s(\Gw).$$ 
This implies \eqref{intL} and therefore $u$ is the unique solution of \eqref{Lp}. Since $s \in (\frac{1}{2},1)$, by \rprop{trG}, $\tr_s(u)=\tr_s(\BBM_s^\Gw[\mu])=\mu$. Finally, estimate \eqref{esL} follows from \rlemma{G} and \rlemma{M}. \qeda
\section{Nonlinear problems}

In this section, we study the nonlinear problem \eqref{Np}. The definition of  weak solutions of \eqref{Np} is given in \rdef{defNp}. 

\subsection{Subcritical absorption}
\noindent \textbf{Proof of Theorem B.} 

\noindent {\sc Monotonicity.} Let $\gt,\gt' \in L^1(\Gw,\gr^s)$, $\gm,\gm' \in L^1(\prt \Gw)$ and $u$ and $u'$ be the weak solutions of \eqref{Np} with data $(\gt,\gm)$ and $(\gt',\gm')$ respectively. We will show that if $\gt \leq \gt'$ and $\gm \leq \gm'$ then $u \leq u'$ in $\Gw$. Indeed, put $v:=(u-u')^+$, it is sufficient to prove that $v \equiv 0$. Since \eqref{form} holds, it follows
$$ |u| \leq \BBG_s^\Gw[|\gt|+|f(u)|] + \BBM_s^\Gw[|\gm|] \quad \text{in } \Gw.  $$
Similarly
$$ |u'| \leq \BBG_s^\Gw[|\gt'|+|f(u')|] + \BBM_s^\Gw[|\gm'|] \quad \text{in } \Gw.  $$
Therefore
$$ 0 \leq v \leq |u| + |u'| \leq \BBG_s^\Gw[|\gt| + |\gt'| + |f(u)| + |f(u')|]  + \BBM_s^\Gw[|\gm|+|\gm'|]:=w.
$$
By Kato inequality, the assumption $\gt \leq \gt'$ and the monotonicity of $f$, we obtain
$$ \fw v \leq \sign^+(u-u')(\gt - \gt') - \sign^+(u-u')(f(u)-f(u')) \leq 0.$$
Therefore
$$ \fw v \leq 0 \leq \fw w \quad \text{in } \Gw. $$
Since $\gm \leq \gm'$, it follows that $\tr_s(v)=0$. By \rlemma{sub}, $v=0$ and thus $u \leq u'$. \medskip

\noindent {\sc Existence.}

\textit{Step 1: Assume that $\gt \in L^\infty(\Gw)$ and  $\mu \in L^\infty(\prt \Gw)$}. 

Put $\hat f(t):=f(t+\BBM_s^\Gw[\mu])-f(\BBM_s^\Gw[\mu])$ and $\hat \gt:=\gt - f(\BBM_s^\Gw[\mu])$. Then $\hat f$ is nondecreasing and $t \hat f(t) \geq 0$ for every $t \in \BBR$ and $\hat \gt \in L^1(\Gw,\gr^s)$. Consider the problem
\bel{ch1} \left\{ \BAL \fw v + \hat f (v) &= \hat \gt \qquad &&\text{in } \Gw \\ 
v &= 0 \qquad &&\text{in } \Gw^c.
\EAL \right. \ee 
By \cite[Proposition 3.1]{CV0} there exists a unique weak solution $v$ of \eqref{ch1}. It means that $v \in L^1(\Gw)$, $\hat f(v) \in L^1(\Gw,\gr^s)$ and 
\bel{ch2} \int_{\Gw} (v \fw \xi + \hat f(v) \xi)\,dx = \int_{\Gw}\xi \hat \gt dx,  \forevery \xi \in \BBX_s(\Gw). \ee
Put $u:=v+\BBM_s^\Gw[\mu]$ then $u \in L^1(\Gw)$ and $f(u) \in L^1(\Gw,\gr^s)$. By \eqref{ch2}  $u$ satisfies \eqref{intN}. 

\medskip

\textit{Step 2: Assume that $0 \leq \gt \in L^1(\Gw,\gr^s)$ and $0 \leq \mu \in L^1(\prt \Gw)$}. 

Let $\{ \gt_n \} \sbs C^1(\ovl \Gw)$ be a nondecreasing sequence converging to $\gt$ in $L^1(\Gw,\gr^s)$ and  $\{ \mu_n \} \sbs C^1(\prt \Gw)$ be a nondecreasing sequence converging to $\mu$ in $L^1(\prt \Gw)$. Then $\{ \BBM_s^\Gw[\mu_n]\}$ is nondecreasing and by \rlemma{M} (ii) it converges to $\BBM_s^\Gw[\mu]$ a.e. in $\Gw$ and in $L^p(\Gw,\gr^s)$ for every $1 \leq p < p_2^*$. Let $u_n$ be the unique solution of \eqref{Np} with $\gt$ and $\gm$ replaced by $\gt_n$ and $\gm_n$ respectively. By step 1 and the monotonicity of $f$, we derive that $\{ u_n \}$ and $\{ f(u_n) \}$ are nondecreasing. Moreover
\bel{en}   \int_{\Gw} (u_n \fw \xi + f(u_n) \xi)\,dx = \int_{\Gw}\xi d\gt_n + \int_{\Gw} \BBM_s^\Gw[\gm_n] \fw \xi \,dx  \qquad \forall \xi \in \BBX_s(\Gw).
\ee
Let $\eta \in C(\ovl \Gw)$ be the solution of
\bel{eta} \left\{ \BAL \fw \eta &= 1 \qquad \text{in } \Gw \\ 
\eta &= 0 \qquad \text{in } \Gw^c,
\EAL \right. \ee
then $c^{-1}\gr^{s}< \eta < c\gr^s$ in $\Gw$ for some $c>1$.
By choosing $\xi=\eta$ in \eqref{en}, we get
\bel{estg} \BAL \norm{u_n}_{L^1(\Gw)} + \norm{f(u_n)}_{L^1(\Gw,\gr^s)} &\leq c(\norm{\gt_n}_{L^1(\Gw,\gr^s)} + \norm{\gm_n}_{L^1(\prt \Gw)}) \\
& \leq c'(\norm{\gt}_{L^1(\Gw,\gr^s)}+ \norm{\gm}_{L^1(\prt \Gw)}).
\EAL \ee
Hence $\{u_n \}$ and $\{ f(u_n) \}$ are uniformly bounded in $L^1(\Gw)$ and $L^1(\Gw,\gr^s)$ respectively. By the monotone convergence theorem, there exists $u \in L^1(\Gw)$ such that $u_n \to u$ in $L^1(\Gw)$ and $f(u_n) \to f(u)$ in $L^1(\Gw,\gr^s)$. By letting $n \to \infty$ in \eqref{en}, we deduce that $u$ satisfies \eqref{intN}, namely $u$ is a weak solution of \eqref{Np}. 

The uniqueness follows from the monotonicity.

\medskip

\textit{Step 3: Assume that $\gt \in L^1(\Gw,\gr^s)$ and $\mu \in L^1(\prt \Gw)$}.

Let $\{ \gt_n \} \sbs C^1(\ovl \Gw)$ be a sequence such that $\{ \gt_n^+ \}$ and $\{ \gt_n^- \}$ are nondecreasing and $\gt_n^{\pm} \to \gt^{\pm}$ in $L^1(\Gw,\gr^s)$. Let $\{ \gm_n \} \sbs C^1(\prt \Gw)$ be a sequence such that $\{ \gm_n^+ \}$ and $\{ \gm_n^- \}$ are nondecreasing and $\gm_n^{\pm} \to \gm^{\pm}$ in $L^1(\prt \Gw)$. Let $u_n$ be the unique weak solution of \eqref{Np} with data $(\gt_n,\gm_n)$, then
\bel{formun}u_n = \BBG_s^\Gw[\gt_n-f(u_n)] + \BBM_s^\Gw[\mu_n]. \ee
Let $w_{1,n}$ and $w_{2,n}$ be the unique weak solutions of \eqref{Np} with data $(\gt_n^+,\gm_n^+)$ and $(-\gt_n^-,-\gm_n^-)$ respectively. Then
\bel{estwi} \norm{w_{i,n}}_{L^1(\Gw)} + \norm{f(w_{i,n})}_{L^1(\Gw,\gr^s)} 
\leq c'(\norm{\gt}_{L^1(\Gw,\gr^s)}+ \norm{\gm}_{L^1(\prt \Gw)}), \quad i=1,2.
\ee
Moreover, for any $n \in \BBN$, $w_{2,n} \leq 0 \leq w_{1,n}$ and
\bel{un2side} -\BBG_s^\Gw[\gt_n^-] - \BBM_s^\Gw[\gm_n^-] \leq w_{2,n} \leq u_n \leq w_{1,n} \leq  \BBG_s^\Gw[\gt_n^+] + \BBM_s^\Gw[\gm_n^+]. \ee
It follows that 
\bel{domi} |u_n| \leq w_{1,n}-w_{2,n} \quad \text{and} |f(u_n)| \leq f(w_{1,n})-f(w_{2,n}). \ee
This, together with \eqref{estwi}, implies
\bel{estun1} \norm{u_{n}}_{L^1(\Gw)} + \norm{f(u_{n})}_{L^1(\Gw,\gr^s)} 
\leq c''(\norm{\gt}_{L^1(\Gw,\gr^s)}+ \norm{\gm}_{L^1(\prt \Gw)}).
\ee
Put $v_n:=\BBG_s^\Gw[\gt_n-f(u_n)]$. By \eqref{estun1}, the sequence $\{ \gt_n-f(u_n \}$ is uniformly bounded in $L^1(\Gw,\gr^s)$. Hence by \cite[Proposition 2.6]{CV1}, the sequence $\{ v_n \}$ is relatively compact in $L^q(\Gw)$ for $1 \leq q < \frac{N}{N-s}$. Consequently, up to a subsequence, $\{ v_n \}$ converges in $L^q(\Gw)$ and a.e. in $\Gw$ to a function $v$. On the other hand, by  \rlemma{M} ii), up to a subsequence, $\{ \BBM_s^\Gw[\mu_n] \}$ converges in $L^q(\Gw,\gr^s)$ for $1\leq q < p_2^*$ and a.e. in $\Gw$ to $\BBM_s^\Gw[\mu]$. Due to \eqref{formun}, we deduce that $\{u_n\}$ converges a.e. in $\Gw$ to $u=v+\BBM_s^\Gw[\mu]$. Since $f$ is continuous, $\{f(u_n)\}$ converges a.e. in $\Gw$ to $f(u)$.

By step 2, the sequences $\{ w_{1,n} \}$, $\{ f(w_{1,n}) \}$, $\{ -w_{2,n} \}$ and $\{ - f(w_{2,n}) \}$ are increasing and converge to $w_1$ in $L^1(\Gw)$, $f(w_1)$ in $L^1(\Gw,\gr^s)$, $-w_2$ in $L^1(\Gw)$ and $-f(w_2)$ in $L^1(\Gw,\gr^s)$ respectively. In light of \eqref{domi} and the generalized dominated convergence theorem, we obtain that $\{ u_n \}$ and $\{ f(u_n) \}$ converge to $u$ and $f(u)$ in $L^1(\Gw)$ and $L^1(\Gw,\gr^s)$ respectively. By passing to the limit in \eqref{en}, we derive that $u$ satisfies \eqref{intN}. \smallskip

The uniqueness follows from the monotonicity. \qeda \medskip

Define 
$$ C(\ovl \Gw, \gr^{-s}):=\{ \zeta \in C(\ovl \Gw): \gr^{-s}\zeta \in C(\ovl \Gw)\}. $$
This space is endowed with the norm
$$ \norm{\zeta}_{C(\ovl \Gw,\gr^{-s})} = \norm{\gr^{-s}\zeta}_{C(\ovl \Gw)}. $$
We say that a sequence $\{ \gt_n \} \sbs \GTM(\Gw,\gr^s)$ converges weakly to a measure $\gt \in \GTM(\Gw,\gr^s)$ if 
$$ \lim_{n \to \infty}\int_{\Gw}\zeta d\gt_n = \int_{\Gw}\zeta d\gt \forevery \zeta \in   C(\ovl \Gw, \gr^{-s}). $$

\noindent \textbf{Proof of Theorem C}.

\noindent {\sc Monotonicity.} The monotonicity can be proved by using a similar argument as in the proof of Theorem B. \smallskip

\noindent {\sc Existence.}
Let $\{ \gt_n\} \sbs C^1(\Gw)$ and $\{\gm_n\} \sbs C^1(\prt \Gw)$ such that $\gt_n^\pm \to \gt^\pm$ weakly  and $\gm_n^\pm \to \gm^\pm$ weakly.  Then there is a positive constant $c$ independent of $n$ such that 
\bel{upseq} \norm{\gt_n}_{\GTM(\Gw,\gr^s)} \leq c\norm{\gt}_{\GTM(\Gw,\gr^s)} \quad \text{and} \quad \norm{\gm_n}_{\GTM(\prt \Gw)} \leq c\norm{\gm}_{\GTM(\prt \Gw)}.
\ee
Let $u_n$,  $w_{1,n}$ and $w_{2,n}$ as in the proof of Theorem B. Then
\bel{un1} |u_n| \leq \max(w_{1,n},-w_{2,n}) \leq \BBG_s^\Gw[|\gt_n|] + \BBM_s^\Gw[|\gm_n|]. \ee
This, together with \eqref{estG1}, \eqref{estM} and \eqref{upseq}, implies that
\bel{un2} \BAL \norm{u_n}_{M^{p_2^*}(\Gw,\gr^s)} \leq c(\norm{\gt_n}_{\GTM(\Gw,\gr^s)}+ \norm{\gm_n}_{\GTM(\prt \Gw)}) 
\leq c'(\norm{\gt}_{\GTM(\Gw,\gr^s)}+ \norm{\gm}_{\GTM(\prt \Gw)}).
\EAL \ee
We have
\bel{intNn} \BAL &\int_{\Gw} (w_{1,n} \fw \xi + f(w_{1,n}) \xi)\,dx = \int_{\Gw}\xi d\gt_n^+ + \int_{\Gw} \BBM_s^\Gw[\gm_n^+] \fw \xi \,dx, \\
&\int_{\Gw} (w_{2,n} \fw \xi + f(w_{2,n}) \xi)\,dx = -\int_{\Gw}\xi d\gt_n^- - \int_{\Gw} \BBM_s^\Gw[\gm_n^-] \fw \xi \,dx\, \forevery \xi \in \BBX_s(\Gw). 
\EAL \ee
From this it follows
\bel{intNn1}  \int_{\Gw} [(w_{1,n} - w_{2,n}) + (f(w_{1,n})-f(w_{2,n}) \eta]\,dx = \int_{\Gw}\eta \, d|\gt_n| + \int_{\Gw} \BBM_s^\Gw[|\gm_n|]\,dx. 
\ee
We infer from \eqref{domi} and the estimate $c^{-1}\gr^s \leq \eta \leq c\gr^s$ that
\bel{estg} \BAL \norm{u_n}_{L^1(\Gw)} + \norm{f(u_n)}_{L^1(\Gw,\gr^s)} &\leq c(\norm{\gt_n}_{L^1(\Gw,\gr^s)} + \norm{\gm_n}_{\GTM(\prt \Gw)}) \\
& \leq c'(\norm{\gt}_{\GTM(\Gw,\gr^s)}+ \norm{\gm}_{\GTM(\prt \Gw)}).
\EAL \ee
This implies that $\{u_n\}$ and $\{ f(u_n) \}$ are uniformly bounded in $L^1(\Gw)$ and $L^1(\Gw,\gr^s)$ respectively. By a similar argument as in step 3 of the proof of Theorem B, we deduce that, up to a subsequence, $\{u_n\}$ converges a.e. in $\Gw$ to a function $u$ and $\{f(u_n)\}$ converges a.e. in $\Gw$ to $f(u)$. By H\"older inequality, we infer that $\{u_n\}$ is uniformly integrable in $L^1(\Gw)$. 

Next we prove that $\{f \circ u_n\}$ is uniformly integrable in $L^1(\Gw,\gr^s)$. Define $\tl f(s):=f(|s|)-f(-|s|)$, $s \in \BBR$. Then $\tl f$ is nondecreasing in $\BBR$ and $|f(s)| \leq \tl f(s)$ for every $s \in \BBR$. For $\ell>0$ and $n \in \BBN$, set
$$ A_n(\ell):=\{ x \in \Gw: |u_n(x)|>\ell  \}, \qquad a_n(\ell):=\int_{A_n(\ell)}\gr^s \,dx. $$  
We take an arbitrary Borel set  $D \sbs \Gw$ and estimate
\bel{esE1} \BAL \int_{D}|f(u_n)|\gr^s dx &=  \int_{D \cap A_n(\ell)}|f(u_n)|\gr^s dx + \int_{D \sms A_n(\ell)}|f(u_n)|\gr^s dx \\
&\leq \int_{A_n(\ell)}\tl f(u_n)\gr^s dx + \tl f(\ell) \int_D \gr^s dx.
\EAL \ee 
On one hand, we have
$$ \int_{A_n(\ell)}\tl f(u_n)\gr^s dx = a_n(\ell)\tl f(\ell) + \int_{\ell}^\infty a_n(s)d \tl f(s). $$
From \eqref{un2}, we infer $a_n(s) \leq \tl c \,s^{-p_2^*}$ where $\tl c$ is a positive constant independent of $n$. Hence, for any $l>\ell$,
\bel{l} \BAL a_n(\ell)\tl f(\ell) + \int_{\ell}^l a_n(s)d \tl f(s) &\leq  \tl c \, \ell^{-p_2^*}\tl f(\ell) + \tl c \int_{\ell}^{l} s^{-p_2^*}d\tl f(s) \\
&\leq \tl c\, l^{-p_2^*} \tl f(l) + \frac{\tl c}{p_2^* +1}\int_{\ell}^l s^{-1-p_2^*}\tl f(s) ds.
\EAL \ee
By assumption \eqref{subcri}, there exists a sequence $\{ l_k \}$ such that $l_k \to \infty$  and $l_k^{-p_2^*}\tl f(l_k) \to 0$ as $k \to \infty$. Taking $l=l_k$ in \eqref{l} and then letting $k \to \infty$, we obtain
\bel{lin}  a_n(\ell)\tl f(\ell) + \int_{\ell}^\infty a_n(s)d \tl f(s) \leq \frac{\tl c}{p_2^* +1}\int_{\ell}^\infty s^{-1-p_2^*}\tl f(s) ds.
\ee
From assumption \eqref{subcri}, we see that the right hand-side of \eqref{lin} tends to $0$ as $\ell \to \infty$. Therefore, for any $\ge>0$, one can choose $\ell>0$ such that the right hand-side of \eqref{lin} is smaller than $\ge/2$. Fix such $\ell$, one then can choose $\gd>0$ small such that if $\int_{D}\gr^s dx < \gd$ then $\tl f(\ell) \int_{D} \gr^s dx < \ge/2$. Therefore, from \eqref{esE1}, we derive
$$ \int_{D} \gr^s dx < \gd  \Lra \int_{D}|f(u_n)|\gr^s dx < \ge. $$
This means $\{f \circ u_n\}$ is uniformly integrable in $L^1(\Gw,\gr^s)$.

By Vitali convergence theorem, we deduce that, up to a subsequence, $u_n \to u$ in $L^1(\Gw)$ and $f(u_n) \to f(u)$ in $L^1(\Gw,\gr^s)$. Since $u_n$ satisfies \eqref{en}, by passing to the limit, we deduce that $u$ is a weak solution of \eqref{Np}. \smallskip

\noindent {\sc Stability.} Assume $\{\gt_n\} \sbs \GTM(\Gw,\gr^s)$ converges weakly to $\gt \in \GTM(\Gw,\gr^s)$ and $\{\gm_n\} \sbs \GTM(\prt \Gw)$ converges weakly to $\gm \in \GTM(\prt \Gw)$. Let $u$ and $u_n$ be the unique weak solution of \eqref{Np} with data $(\gt,\gm)$ and $(\gt_n,\gm_n)$ respectively. Then by a similar argument as in Existence part, we deduce that $u_n \to u$ in $L^1(\Gw)$ and $f(u_n) \to f(u)$ in $L^p(\Gw,\gr^s)$. \qeda.

\bprop{GfM} Assume $f$ is a continuous nondecreasing function on $\BBR$ satisfying $f(0)=0$ and \eqref{subcri}. Then for every $z \in \prt \Gw$,
\bel{limGfM} \lim_{\Gw \ni x \to z}\frac{\BBG_s^\Gw[f(M_s^\Gw(\cdot,z))](x)}{M_s^\Gw(x,z)}=0.
\ee
\es
\proof By \eqref{Ga}, 
$$ G_s^\Gw(x,y) \leq c_{14} \gr(x)^s |x-y|^{-N}\min\{ \gr(y)^s, |x-y|^s \}\qquad \forall x\neq y. $$
Hence 
\bel{GfM1} \BAL \frac{\BBG_s^\Gw[f(M_s^\Gw(\cdot,z))](x)}{M_s^\Gw(x,z)} \leq c_{15}|x-z|^N \int_{\Gw}|x-y|^{-N}\min\{ |x-y|^s, |y-z|^s \}f(|y-z|^{s-N})dy.
\EAL \ee
Put
\bel{CD} \CD_1:=\Gw \cap B(x,|x-z|/2), \quad \CD_2:=\Gw \cap B(z,|x-z|/2), \quad \CD_3:=\Gw \sms (\CD_1 \cup \CD_2), \ee
$$ I_i:=|x-z|^N \int_{\CD_i}|x-y|^{-N}\min\{ |x-y|^s, |y-z|^s \}f(|y-z|^{s-N})dy, \quad i=1,2,3.
$$
For every $y \in \CD_1$, $|x-z| \leq 2|y-z|$, therefore
$$ I_1 \leq c_{16}|x-z|^Nf(|x-z|^{s-N})\int_{\CD_1}|x-y|^{s-N}dy \leq c_{17}|x-z|^{N+s}f(|x-z|^{s-N}). $$
Hence 
\bel{I11} \lim_{x \to z}I_1 \leq c_{17} \lim_{x \to z}|x-z|^{N+s}f(|x-z|^{s-N}) = 0. \ee
We next estimate $I_2$. For every $y \in \CD_2$, $|x-z| \leq 2|x-y|$, hence
$$ I_2 \leq c_{18} \int_{\CD_2}|y-z|^s f(|y-z|^{s-N})dy \leq c_{37} \int_{|x-z|^{s-N}}^\infty t^{-1-p_2^*}f(t)dt.
$$
Therefore, by \eqref{subcri},
\bel{I21} \lim_{x \to z}I_2 \leq c_{19} \lim_{x \to z}\int_{|x-z|^{s-N}}^\infty t^{-1-p_2^*}f(s)ds = 0.
\ee
Finally, we estimate $I_3$. For every $y \in \CD_3$, $|y-z| \leq 3|x-y|$, therefore
\bel{I31} I_3 \leq c_{20}|x-z|^N \int_{\CD_3} |y-z|^{s-N} f(|y-z|^{s-N}) dy \leq c_{21} |x-z|^N \int_0^{|x-z|^{s-N}} t^{-\frac{N}{N-s}}f(t)dt.
\ee
Put 
$$g_1(r)=\int_0^{r^{s-N}} t^{-\frac{N}{N-s}}f(t)dt, \quad g_2(r)=r^{-N}. $$
If $\lim_{r \to 0}g_1(r)< \infty$, then $\lim_{x \to z}I_3=0$ by \eqref{I31}. Otherwise, $\lim_{r \to 0}g_1(r)=\infty=\lim_{r \to 0}g_2(r)$. Therefore, by L' H\^opital's rule,
\bel{I32} \lim_{r \to 0}\frac{g_1(r)}{g_2(r)}=\lim_{r \to 0}\frac{g'_1(r)}{g'_2(r)}=\lim_{r \to 0}\frac{N-s}{N}r^{N+s}f(r^{s-N}) = 0. \ee
By combining \eqref{I31} and \eqref{I32} we obtain
\bel{I33} \lim_{x \to z}I_3 \leq c_{22} \lim_{x \to z}|x-z|^N \int_0^{|x-z|^{s-N}} t^{-\frac{N}{N-s}}f(t)dt = 0. \ee
We deduce \eqref{GfM1} by gathering \eqref{I11}, \eqref{I21} and \eqref{I33}. \qeda\smallskip

\noindent \textbf{Proof of Theorem D}. From Theorem C we get
\bel{beh} kM_s^\Gw(x,z) - \BBG_s^\Gw[f(M_s^\Gw(\cdot,z))](x) \leq u_{k,z}^\Gw(x) \leq kM_s^\Gw(x,z), \ee
which implies
$$ k - \frac{\BBG_s^\Gw[f(M_s^\Gw(\cdot,z))](x)}{M_s^\Gw(x,z)} \leq \frac{u_{k,z}^\Gw(x)}{M_s^\Gw(x,z)} \leq k.$$
We derive \eqref{ukz} due to \rprop{GfM}. \qeda
\subsection{Power absorption}
In this subsection we assume that $0 \in \prt \Gw$. Let $0<p<p_2^*$ and denote by $u_k^\Gw$ the unique solution of \eqref{Nppd}. By Theorem C, $u_k^\Gw \leq kM_s^\Gw(\cdot,0)$ and  $k \mapsto u_k^\Gw$ is increasing. Therefore, it is natural to investigate $\lim_{k \to \infty}u_k^\Gw$.

For any $\ell>0$, put 
$$ T_\ell[u](y):=\ell^{\frac{2s}{p-1}}u(\ell y), \quad  y \in \Gw_\ell:=\ell^{-1}\Gw.$$ 
If $u$ is a solution of \eqref{eqa} in $\Gw$ then $T_\ell[u]$ is a solution of \eqref{eqa} in $\Gw_\ell$. 

By \rcor{rad}, the function 
\bel{U} x\mapsto U(x)=\ell_{s,p}|x|^{-\frac{2s}{p-1}}, \quad x \neq 0, \ee 
where $\ell_{s,p}$ is a positive constant, is a radial singular solution of 
\bel{Nep} \fw u + u^p = 0 \quad \text{in } \BBR^N\setminus\{0\}. \ee

\blemma{upbound} Assume $p \in (p_1^*,p_2^*)$. Then there exists a positive constant $C$ depending on $N$, $s$, $p$ and the $C^2$ characteristic of $\Gw$ such that the following holds. If $u$ is a positive solution of \eqref{eqa} satisfying $u \leq U$ in $\Gw$ then there holds
\bel{upb1} u(x) \leq C\gr(x)^s|x|^{-\frac{(p+1)s}{p-1}} \quad \forall x \in \Gw.
\ee
\es
\proof Let $P \in (\prt \Gw \sms \{0\}) \cap B_1(0)$ and put
$$ d=d(P):=\frac{1}{2}|P|< \frac{1}{2}. $$
Put
$$ u_d(y)=T_d[u](y), \quad y \in \Gw_d:=d^{-1}\Gw.
$$
Then $u_d$ is a solution of 
\bel{eqOd} \left\{ \BA{rll} \fw u + u^p &= 0 \qquad &\text{in } \Gw_d \\ 
u &= 0 \qquad &\text{in } (\Gw_d)^c.
\EA \right. \ee
Moreover 
$$ u_d(y) \leq T_d[U](y)=d^{\frac{2s}{p-1}}U(dy)=\ell_{s,p}|y|^{-\frac{2s}{p-1}}=U(y). $$
Put $P_d=d^{-1}P$ and let $\gb_0$ be the constant in \rprop{beta0}. We may assume $\gb_0 \leq \frac{1}{4}$. Let $\zeta_P \in C^\infty(\BBR^N)$ such that $0 \leq \zeta \leq 1$ in $\BBR^N$, $\zeta=0$ in $B_{\gb_0}(P_d)$ and $\zeta=1$ in $\BBR^N \sms B_{2\gb_0}(P_d)$. Let $\eta_d \in C(\ovl \Gw_d)$ be the solution of \eqref{eta} with $\Gw$ replaced by $\Gw_d$. For $l>0$, denote
$$ V_{d,l}:= \zeta_P \, U + l \, \eta_d. $$  
We will compare $u_d$ with $V_{d,l}$. \smallskip

\noindent \textit{Step 1: We show that $V_{d,l}$ is a super solution of \eqref{eqOd} for $l$ large enough.}

For $y \in \Gw_d \sms B_{4\gb_0}(P_d)$, $\zeta_P(y)=1$ and hence
$$ \BAL
\fw (\zeta_P U)(y) &= \lim_{\ge \to 0}\int_{\BBR^N \sms B_\ge(y)}\frac{U(y)-\zeta_P(z)U(z)}{|y-z|^{N+2s}}dz \\
&=\fw U(y)  + \lim_{\ge \to 0}\int_{\BBR^N \sms B_\ge(y)}\frac{U(z)-\zeta_P(z)U(z)}{|y-z|^{N+2s}}dz \\
&\geq \fw U(y) - \int_{B_{\frac{1}{2}}(P_d)}\frac{U(z)}{|y-z|^{N+2s}}dz \\
&\geq \fw U(y) - c_{26},
\EAL $$
where $c_{26}=c_{26}(N,s,p,\gb_0)$.
Since $(\Gw_d \cap B_{2\gb_0}(0)) \sbs (\Gw_d \sms B_{4\gb_0}(P_d))$, it follows that, for any $y \in \Gw_d \cap B_{2\gb_0}(0) \sms \{0\}$,
$$ \BAL
\fw V_{d,l}(y) + (V_{d,l}(y))^p &=\fw (\zeta_P U)(y) + l \fw \eta_d(y) + (\zeta_P(y) U(y) + l \eta_d(y))^p \\
&\geq \fw U(y) - c_{26} + l + U(y)^p. 
\EAL $$
Therefore if we choose $l \geq c_{26}$ then 
\bel{bigV1}  \fw V_{d,l} + (V_{d,l})^p \geq 0 \quad \text{in } \Gw_d \cap B_{2\gb_0}(0) \sms \{0\}. \ee
Next we see that there exists $c_{27}>0$ such that
$$ |\fw (\zeta_P U)| \leq c_{27} \quad \text{in } \Gw_d \sms B_{2\gb_0}(0). $$
Consequently,
$$ \BAL \fw V_{d,l} &= \fw (\zeta_P U) + l \fw \eta_d \\
&\geq -c_{27} + l.
\EAL $$
Therefore if we choose $l \geq c_{27}$ then 
\bel{bigV2} \fw V_{d,l} \geq 0 \quad \text{in } \Gw_d \sms B_{2\gb_0}(0). \ee
By combining \eqref{bigV1} and \eqref{bigV2}, for $l\geq \max\{ c_{26}, c_{27} \}$, we deduce that $V_{d,l}$ is a super solution of \eqref{eqOd}. \smallskip

\noindent \textit{Step 2: We show that $u_d \leq V_{d,l}$ in $\Gw_d$.} By contradiction, we assume that there exists $x_0 \in \Gw_d$ such that 
$$ (u_d - V_{d,l})(x_0) = \max_{x \in \Gw_d}(u_d - V_{d,l}) > 0. $$
Then $\fw (u_d-V_{d,l})(x_0) \geq 0$. It follows that
$$ 0 \leq \fw (u_d - V_{d,l})(x_0) \leq  -(u_d(x_0)^p - V_{d,l}(x_0)^p)<0. $$
This contradiction implies that $u_d \leq V_{d,l}$ in $\Gw_d$. \smallskip

\noindent \textit{Step 3: End of proof.} From step 2, we deduce that
$$ u_d \leq l\eta_d \quad \text{in } \Gw_d \cap B_{\gb_0}(P_d). $$
We note that $\eta_d(y) \leq c\dist(y,\prt \Gw_d)^s$ for every $y \in \Gw_d$. Here the constant $c$ depends on $N$, $s$ and the $C^2$ characteristic of $\Gw_d$. Since $d < \frac{1}{2}$, a $C^2$ characteristic of $\Gw_d$ can be taken as a $C^2$ characteristic of $\Gw$. Therefore the constant $c$ can be taken independently of $P$. Consequently,
$$ u_d(y) \leq l c \,\dist(y, \prt \Gw_d)^s \quad \forall y \in \Gw_d \cap B_{\gb_0}(P_d). $$
This implies
\bel{upp1} u(x) \leq c' \gr(x)^s d^{-\frac{(p+1)s}{p-1}} \quad \forall x \in \Gw \cap B_{d\gb_0}(P).
\ee
Put 
$$ \CF_1:=\Gw_{\gb_0} \cap B_{\frac{1}{1+\gb_0}}(0) \cap \{x: \gr(x) \leq \gb_0 |x| \}, \quad \CF_2:=\Gw_{\gb_0} \cap B_{\frac{1}{1+\gb_0}}(0) \cap \{x: \gr(x) > \gb_0 |x| \}.
$$ 
If $x \in \CF_1$ then let $P \in \prt \Gw \sms \{0\}$ such that $\gr(x)=|x-P|$. It follows that
\bel{upp2} \frac{1}{2}(1-\gb_0)|x|<d=\frac{1}{2}|P|  \leq \frac{1}{2}(1+\gb_0)|x|<\frac{1}{2}. \ee
By combining \eqref{upp1} and \eqref{upp2}, we get
$$ u(x) \leq c'(1-\gb_0)^{-\frac{(p+1)s}{p-1}} \gr(x)^s |x|^{-\frac{(p+1)s}{p-1}}. $$
If $x \in \CF_2$ then \eqref{upb1} follows from the assumption $u \leq U$. Thus \eqref{upb1} holds for every $x \in \Gw_{\gb_0} \cap B_{\frac{1}{1+\gb_0}}(0)$. If $x \in \Gw \sms B_{\frac{1}{1+\gb_0}}(0)$ then by a similar argument as in Step 1 and Step 2 without similarity transformation, we deduce that there exist constants $c$ and  $\tl \gb \in (0,\frac{1}{2(1+\gb_0)})$ depending on $N$, $s$, $p$ and the $C^2$ characteristic of $\Gw$ such that \eqref{upb1} holds in $B_{\tl \gb}(P) \cap \Gw$ for every $P \in \prt \Gw \sms B_{\frac{1}{1+\gb_0}}(0)$. Finally, since $u \leq U$, \eqref{upb1} holds in $D_{\frac{\tl \gb}{2}}=\{x \in \Gw: \gr(x) > \frac{\tl \gb}{2}\}$. Thus \eqref{upb1} holds in $\Gw$. \qeda 

\blemma{GMp} Let $p \in (0,p_2^*)$. There exists a constant $c=c(N,s,p,\Gw)>0$ such that for any $x \in \Gw$ and $z \in \prt \Gw$, there holds
\bel{GMp} \BBG_s^\Gw[M_s^\Gw(\cdot,z)^p](x) \leq \left\{ \BAL &c\gr(x)^s |x-z|^{s-(N-s)p} \qquad &&\text{if} \quad \frac{s}{N-s}<p<p_2^* \\
-&c\gr(x)^s\ln|x-z| &&\text{if} \quad p=\frac{s}{N-s} \\
&c\gr(x)^s &&\text{if} \quad 0<p<\frac{s}{N-s}.
\EAL \right.
\ee
\es
\proof We use a similar argument as in the proof of \rprop{GfM}. It is easy to see that for every $x \in \Gw$ and $z \in \prt \Gw$,
\bel{Mq1}
\BBG_s^\Gw[M_s^\Gw(\cdot,z)^p](x) \leq c_{23}\gr(x)^s\int_{\Gw}|x-y|^{-N}|y-z|^{(s-N)p}\min\{ |x-y|^s, |y-z|^s \}dy
\ee
Let $\CD_i$, $i=1,2,3$ be as in \eqref{CD} and put
$$ J_i:=\gr(x)^s\int_{\CD_i}|x-y|^{-N}|y-z|^{(s-N)p}\min\{ |x-y|^s, |y-z|^s \}dy. $$
By proceeding as in the proof of \rprop{GfM} we deduce easily that there is  positive constants $c_{24}=c_{24}(N,s,p,\Gw)$ such that
\bel{J12} J_i \leq c_{24}\gr(x)^s |x-z|^{s-(N-s)p}, \quad i=1,2, \ee
and 
\bel{J3} J_3 \leq c_{24}\gr(x)^s\int_{|x-z|/2}^{diam(\Gw)}r^{s-1-(N-s)p}dr \leq \left\{ \BAL
&c_{25}\gr(x)^s |x-z|^{s-(N-s)p} &&\text{if} \quad \frac{s}{N-s}<p<p_2^* \\
-&c_{25}\gr(x)^s\ln|x-z| &&\text{if} \quad p=\frac{s}{N-s} \\
&c_{25}\gr(x)^s &&\text{if} \quad 0<p<\frac{s}{N-s}.
\EAL \right. \ee
Combining \eqref{J12} and \eqref{J3} implies \eqref{GMp}. \qeda \medskip

\bprop{uinf1} Assume $p \in (p_1^*,p_2^*)$. Then $u_\infty^\Gw:=\lim_{k \to 0}u_k^\Gw$ is a positive solution of \eqref{eqa}. Moreover, there exists $c=c(N,s,p,\Gw)>0$ such that
\bel{uinf2} c^{-1}\gr(x)^s|x|^{-\frac{(p+1)s}{p-1}} \leq u_\infty^\Gw(x) \leq c\gr(x)^s |x|^{-\frac{(p+1)s}{p-1}} \forevery x \in \Gw.
\ee
\es
\proof We first claim that for any $k>0$,
\bel{ukU} u_k^\Gw \leq U \quad \text{in }  \Gw. \ee 
Indeed, by \eqref{estMar}, 
$$u_k^\Gw(x) \leq kM_s^\Gw(x,0) \leq c_{28}k\gr(x)^s|x|^{-N} \leq c_{28}k |x|^{s-N}  \quad \forall x \in \Gw. $$ 
Since $p <p_2^*$, it follows that 
$$\lim_{\Gw \ni x \to 0}\frac{u_k^\Gw(x)}{U(x)} = 0.$$
By proceeding as in Step 2 of the proof of \rlemma{upbound}, we deduce that $u_k^\Gw \leq U$ in $\Gw$. 

Consequently, $u_\infty^\Gw:=\lim_{k \to \infty}u_k^\Gw$ is a solution of \eqref{eqa} vanishing on $\prt \Gw \sms \{0\}$ and satisfying $u_\infty^\Gw \leq U$ in $\Gw$. In light of \rlemma{upbound}, we obtain the upper bound in \eqref{uinf2}.

Next we prove the lower bound in \eqref{uinf2}. By \eqref{estMar} and \rlemma{GMp}, for any $k >0$ and $x \in \Gw$, we have
$$ \BAL u_k^\Gw(x) &\geq kM_s^\Gw(x,0) - k^p\BBG_s^\Gw[M_s^\Gw(\cdot,0)^p](x) \\
&\geq c_{29}^{-1} k \gr(x)^s|x|^{-N}(1-c_{29}c_{30} k^{p-1}|x|^{N+s- (N-s)p}). \\ 
\EAL $$ 
For $x \in \Gw$, one can choose $r>0$ such that $x \in \Gw \cap (B_{2r}(0) \sms B_r(0))$. Choose $k=ar^{-\frac{N+s-(N-s)p}{p-1}}$, where $a>0$ will be made precise later on, then
$$ u_k^\Gw(x) \geq c_{31} a \, \gr(x)^s |x|^{-\frac{(p+1)s}{p-1}}(1-c_{29}c_{30}a^{p-1}). 
$$
By choosing $a=(2c_{29}c_{30})^{-\frac{1}{p-1}}$, we deduce for any $x \in \Gw$ there exists $k>0$ depending on $|x|$ such that
$$ u_k^\Gw(x) \geq c_{32}\gr(x)^s|x|^{-\frac{(p+1)s}{p-1}}. $$
Since $u_\infty^\Gw \geq u_k^\Gw$ in $\Gw$ we obtain the first inequality in \eqref{uinf2}. \qeda

\bprop{ukl} Assume $0<p\leq p_1^*$. There exist $k_0=k_0(N,s,p)$ and $c=c(N,s,p,\Gw)$ such that the following holds. There exists a decreasing sequence of positive numbers $\{r_k\}$ such that $\lim_{k \to \infty}r_k=0$ and for any $k > k_0$,
\bel{ukl1} u_k^\Gw(x) \geq \left\{ \BAL &c\gr(x)^s |x|^{-N-s} \quad &&\text{if } 0<p <p_1^*, \\[4mm]
&c\gr(x)^s |x|^{-N-s}(-\ln|x|)^{-1} \quad &&\text{if } p=p_1^*,
\EAL \right.
, \forevery x \in \Gw \sms B_{r_k}(0). \ee
\es
\proof  For any $\ell>0$, we have 
\bel{ukl2} u_\ell^\Gw(x) \geq \ell M_s^\Gw(x,0) - \ell^{p}\BBG_s^\Gw[M_s^\Gw(\cdot,0)^p](x) \forevery x \in  \Gw. \ee
\noindent \textit{Case 1: $p \in (\frac{s}{N-s},p_1^*)$}. Put $k_1:=(2c_{29}c_{30})^{\frac{s}{N+2s-Np}}$ and take $k>k_1$. For $\ell>0$, put $r_\ell=\ell^{-\frac{1}{s}}$, then $\ell=r_\ell^{-s}$.  Take arbitrarily $x \in \Gw \sms B_{r_k}(0)$ then one can choose $\ell \in (\max(2^{-s}k,k_1),k)$ such that $x \in \Gw \cap (B_{r_\ell}(0) \sms B_{\frac{r_\ell}{2}}(0))$. From \eqref{ukl2}, \eqref{estMar} and \eqref{GMp}, we get
$$ \BAL u_\ell^\Gw(x) &  \geq c_{29}^{-1}\ell \gr(x)^s |x|^{-N}(1-c_{29}c_{30} \ell^{p-1}|x|^{N+s-(N-s)p}) \\
& \geq  c_{29}^{-1}\gr(x)^s |x|^{-N}r_\ell^{-s}(1 - c_{29}c_{30}r_\ell^{N+2s-Np}) \\
& \geq (2c_{29})^{-1} \gr(x)^s |x|^{-N} r_\ell^{-s} \\
& \geq c_{33} \gr(x)^s |x|^{-N-s}.
\EAL $$
Here the first estimate holds since $\frac{N}{N-s}<p<p_2^*$ and the third estimate holds since $p<p_1^*$ and $\ell >k_1$. Since $k > \ell$, we deduce that
\bel{ukl3} u_k^\Gw(x) \geq c_{33}\gr(x)^s |x|^{-N-s}, \forevery x \in \Gw \sms B_{r_k}(0). \ee

\noindent \textit{Case 2: $p=\frac{s}{N-s}$}. Put $k_2=(\frac{2c_{29}c_{30}(1+s)}{s})^{\frac{s}{N-sp}}$ and take $k>k_2$. For $\ell>0$, put $r_\ell=\ell^{-\frac{1}{s}}$, then $\ell=r_\ell^{-s}$.  Take arbitrarily $x \in \Gw \sms B_{r_k}(0)$ then one can choose $\ell \in (\max(2^{-s}k,k_2),k)$ such that $x \in \Gw \cap (B_{r_\ell}(0) \sms B_{\frac{r_\ell}{2}}(0))$. From \eqref{ukl2}, \eqref{estMar} and \rlemma{GMp}, we get
$$ \BAL u_\ell^\Gw(x) &  \geq c_{29}^{-1}\ell \gr(x)^s |x|^{-N}(1 + c_{29}c_{30} \ell^{p-1}|x|^{N}\ln|x|) \\
& \geq  c_{29}^{-1} \gr(x)^s |x|^{-N}r_\ell^{-s}(1 + c_{29}c_{30}r_\ell^{N+s-sp}\ln(\frac{r_\ell}{2})) \\
& \geq (2c_{29})^{-1} \gr(x)^s |x|^{-N} r_\ell^{-s} \\
& \geq c_{33} \gr(x)^s |x|^{-N-s}.
\EAL $$
Here the third estimate holds since $\ell >k_2$ and $N-sp>0$. Therefore \eqref{ukl3} holds. \smallskip

\noindent \textit{Case 3: $p \in (0,\frac{s}{N-s})$}.  Put $k_3=(2c_{29}c_{30})^{\frac{s}{N+s-sp}}$ and take $k>k_3$. For $\ell>0$, put $r_\ell=\ell^{-\frac{1}{s}}$, then $\ell=r_\ell^{-s}$.  Take arbitrarily $x \in \Gw \sms B_{r_k}(0)$ then one can choose $\ell \in (\max(2^{-s}k,k_3),k)$ such that $x \in \Gw \cap (B_{r_\ell}(0) \sms B_{\frac{r_\ell}{2}}(0))$. From \eqref{ukl2}, \eqref{estMar} and \eqref{GMp}, we get
$$ \BAL u_\ell^\Gw(x) &  \geq c_{29}^{-1}\ell \gr(x)^s |x|^{-N}(1 - c_{29}c_{30} \ell^{p-1}|x|^{N}) \\
& \geq  c_{29}^{-1} \gr(x)^s |x|^{-N}r_\ell^{-s}(1 - c_{29}c_{30}r_\ell^{N+s-sp}) \\
& \geq (2c_{29})^{-1}\gr(x)^s |x|^{-N} r_\ell^{-s} \\
& \geq c_{33} \gr(x)^s |x|^{-N-s}.
\EAL $$
Here the third estimate holds since $\ell >k_3$ and $N+s-sp>0$. Therefore \eqref{ukl3} holds. \smallskip

\noindent \textit{Case 4: $p = p_1^*$}. Put $k_4=\exp( (2c_{29}c_{30})^{\frac{s}{N+s-(N-s)p}})$ and take $k>k_4$. For $\ell>0$, put $r_\ell=(\ell \ln(\ell))^{-\frac{1}{s}}$, then $\ell\ln(\ell)=r_\ell^{-s}$ and $\ell <r_\ell^{-s}$ when $\ell>3$.  Take arbitrarily $x \in \Gw \sms B_{r_k}(0)$ then one can choose $\ell \in (\max(2^{-s}k,k_4),k)$ such that $x \in \Gw \cap (B_{r_\ell}(0) \sms B_{\frac{r_\ell}{2}}(0))$. From \eqref{ukl2}, \eqref{estMar} and \eqref{GMp}, we get
$$ \BAL u_\ell^\Gw(x) &  \geq c_{29}^{-1}\ell \gr(x)^s |x|^{-N}(1-c_{29}c_{30} \ell^{p-1}|x|^{N+s-(N-s)p}) \\
& \geq  c_{29}^{-1} \ell \gr(x)^s |x|^{-N}(1 - c_{29}c_{30} \ell^{p-1} (\ell \ln(\ell))^{-\frac{N+s-(N-s)p}{s}}) \\
&=c_{29}^{-1} \ell \gr(x)^s |x|^{-N}(1 - c_{29}c_{30} \ln(\ell)^{-\frac{N+s-(N-s)p}{s}}) \\
& \geq (2c_{29})^{-1} \ell \gr(x)^s |x|^{-N}  \\
& \geq c_{34} \gr(x)^s |x|^{-N-s}(-\ln|x|)^{-1}.
\EAL $$
Here the inequality holds since $p=p_1^*$ and the last estimate follows from the following estimate
$$ \ell =\frac{r_\ell^{-s}}{\ln(\ell)}> \frac{|x|^{-s}}{-s 2^s \ln|x|}. $$
Since $u_k^\Gw(x) \geq u_\ell^\Gw(x)$, we derive 
$$ u_k^\Gw(x) \geq  c_{34} \gr(x)^s |x|^{-N-s}(-\ln|x|)^{-1}. $$
By putting $k_0:=\max(k_1,k_2,k_3,k_4)$, we obtain \eqref{ukl1}. \qeda \medskip

\bprop{uinf2} Assume $0<p \leq p_1^*$. Then $\lim_{k \to \infty}u_k^\Gw(x)=\infty$ for every $x \in \Gw$.
\es
\proof The proposition can be obtained by adapting the argument in the proof of \cite[Theorem 1.2]{CAHM}. Let $r_0>0$ and put 
$$\gth_k:=\int_{B_{r_0}(0)}u_k^\Gw(x)dx. $$
Then
$$ \gth_k \geq c\int_{(B_{r_0}\cap \Gw) \sms B_{r_k}(0)}\gr(x)^s|x|^{N-s}(-\ln|x|)^{-1}dx,
$$ 
which implies 
\bel{thk} \lim_{k \to \infty}\gth_k=\infty. \ee
Fix $y_0 \in \Gw \sms \ovl B_{r_0}(0)$ and set $\gd:=\frac{1}{2}\min\{ \gr(y_0), |y_0|-r_0  \}$. By \cite[Lemma 2.4]{CY} there exists a unique classical solution $w_k$ of the following problem
\begin{equation}
 \left\{
 \begin{alignedat}{2}
   \fw w_k + w_k^p & = 0 \quad &&\text{in } B_{\delta}(y_0), \\\
   w_k &= 0 &&\text{in } \BBR^N \sms (B_{\delta}(y_0) \cup B_{r_0}(0)), \\
   w_k &= u_k^\Gw &&\text{in } B_{r_0}(0).
 \end{alignedat} 
  \right.
\end{equation}
By \cite[Lemma 2.2]{CY}, 
\bel{wu} u_k^\Gw \geq w_k \quad \text{in } B_{\delta}(y_0).\ee
Next put
$\tl w_k:= w_k - \chi_{B_{r_0}(0)}u_k $
then $\tl w_k = w_k$ in $B_{\delta}(y_0)$. Moreover, for $x \in B_{\delta}(y_0)$
\begin{equation}
\begin{alignedat}{2}
\fw \tl w_k(x) &=  \lim_{\ge \to 0}\int_{B_{\delta}(y_0) \sms B_\ge(x)}\frac{w_k(x)-w_k(z)}{|z-x|^{N+2s}}dz + \lim_{\ge \to 0}\int_{B_{\delta}^c(y_0) \sms B_\ge(x) }\frac{w_k(x)}{|z-x|^{N+2s}}dz \\
&= \lim_{\ge \to 0}\int_{\BBR^N \sms B_\ge(x)}\frac{w_k(x)-w_k(z)}{|z-x|^{N+2s}}dz + \int_{B_{r_0}(0)}\frac{u_k^\Gw(z)}{|z-x|^{N+2s}}dz \\
&\geq \fw w_k(x) + A \theta_k
\end{alignedat} 
\end{equation}
where $A=(|y_0|+r_0)^{-N-2s}$. It follows that, for $x \in B_{\delta}(y_0)$,
\begin{equation}
\fw \tl w_k(x) + \tl w_k^p(x) \geq \fw w_k(x) + w_k^p(x) + A\theta_k = A\theta_k.
\end{equation}
Therefore $\tl w_k \in C(\ovl {B_\gd(y_0)})$ is a supersolution of 
\begin{equation} \label{eAth}
 \left\{
 \begin{alignedat}{2}
   \fw w + w^p & = A\theta_k \quad &&\text{in } B_{\delta}(y_0), \\\
   w &= 0 &&\text{in } \BBR^N \sms B_{\delta}(y_0).  \\
 \end{alignedat} 
  \right.
\end{equation}
Let $\eta_0 \in C(\ovl {B_\gd(y_0)})$ be the unique solution of 
\begin{equation}
 \left\{
 \begin{alignedat}{2}
   \fw \eta_0 & = 1 \quad &&\text{in } B_{\delta}(y_0), \\\
   \eta_0 &= 0 &&\text{in } \BBR^N \sms B_{\delta}(y_0). \\
 \end{alignedat} 
  \right.
\end{equation}
We can choose $k$ large enough so that the function 
$$ \frac{\eta_0(A\theta_k)^{\frac{1}{p}}}{2\max_{\BBR^N}\eta_0} $$
is a subsolution of \eqref{eAth}. By \cite[Lemma 2.2]{CY} we obtain
\begin{equation} \label{wgeq} \tl w_k(x) \geq  \frac{\eta_0(A\theta_k)^{\frac{1}{p}}}{2\max_{\BBR^N}\eta_0} \quad \forall x \in B_{\delta}(y_0).
\end{equation}
Put 
$$ \underline c:=\min_{x \in B_{\delta}(y_0)}\frac{\eta_0}{2\max_{\BBR^N}\eta_0} $$
then we derive from \eqref{wgeq} that 
\bel{wl} w_k(x) \geq \underline c (A\theta_k)^{\frac{1}{p}} \quad \forall x \in B_\gd(y_0).  \ee
By combining \eqref{thk}, \eqref{wu} and \eqref{wl}, we deduce that
$$ \lim_{k \to \infty}u_k^\Gw(x) = \infty \quad \forall x \in B_{\frac{\delta}{2}}(y_0). $$
This implies 
$$ \lim_{k \to \infty}u_k^\Gw(x) = \infty \quad \forall x \in \Gw. $$
\qeda

\bth{unbduk} Assume $p \in (1,p_2^*)$ and  either $\Gw=\BBR_+^N:=\{x=(x',x_N):x_N>0\}$ or $\prt\Gw$ is compact with $0\in\prt\Gw$. Then, for any $k>0$, there exists a unique solution solution $u_k^\Gw$ of problem \eqref{Nppd} satisfying $u_k^\Gw \leq kM_s^\Gw(\cdot,0)$ in $\Gw$ and
\bel{ukO} \lim_{|x| \to 0}\frac{u_k^\Gw(x)}{M_s^\Gw(x,0)}=k. \ee
Moreover, the map $k \mapsto u_k^\Gw$ is increasing.
\es
\proof \textit{Step 1: Existence.} For $R>0$ we set $\Gw_R=\Gw\cap B_R$ and let $u:=u^{\Gw_R}_k$ be the unique solution of 
\bel {ukOR} \left\{ \BAL \fw u + u^p  &= 0 \quad &&\text{in } \Gw_R \\
 		\tr_s(u) &= k\gd_0 &&~ \\
      u &= 0 &&\text{on } \Gw_R^c.
 	\EAL \right. \ee
Then
	 \bel {Pc+2}u^{\Gw_R}_{k}(x)\leq kM_s^{\Gw_R}(x,0)\qq\forall x\in \Gw_R. 
	 \ee	 
Since $R\mapsto M_s^{\Gw_R}(.,0)$ is increasing, it follows from \eqref{ukz} that  $R\mapsto u^{\Gw_R}_{k}$ is increasing too with the limit $u^*$ and there holds 
	\bel {Pc+3}
	u^*(x)\leq kM_s^{\Gw}(x,0)\qq\forall x\in \Gw.
	\ee
From \eqref{Pc+2}, we deduce that 
$$ u_k^{\Gw_R}(x) \leq ck |x|^{s-N} \quad \forevery x \in \Gw_R$$
where $c$ depends only on $N$, $s$ and the $C^2$ characteristic of $\Gw$.  Hence by the regularity up to the boundary \cite{RS}, $\{ u^{\Gw_R}_{k}\}$ is uniformly bounded in $C_{loc}^s(\overline \Gw\setminus B_\ge)$ and in $C_{loc}^{2s+\alpha}(\Gw \sms B_\ge)$ for any $\ge>0$. Therefore, $\{u_k^{\Gw_R}\}$ converges locally uniformly, as $R \to \infty$, to $u^*\in C(\overline\Gw\setminus\{0\}) \cap C^{2s+\alpha}(\Gw)$. Thus $u^*$ is a positive solution of \eqref{eqa}. 
Moreover by combining \eqref{ukz}, \eqref{Pc+2}, the fact that $M_s^{\Gw_R} \uparrow M_s^{\Gw}$ and $u_k^{\Gw_R} \uparrow u_k^{\Gw}$, we deduce that $\tr_s(u^*)=k\gd_0$ and 
$$ \lim_{\Gw \ni x \to 0}\frac{u^*(x)}{M_s^\Gw(x,0)}=k. 
$$
\noindent \textit{Step 2: Uniqueness.} Suppose $u$ and $u'$ are two weak solutions of \eqref{eqa} satisfying $\max\{ u,u'\} \leq kM_s^{\Gw}(\cdot,0)$ in $\Gw$ and
\bel{luba} \lim_{\Gw \ni x\to 0}\frac{u(x)}{M_s^{\Gw}(x,0)}=\lim_{\Gw \ni x\to 0}\frac{u'(x)}{M_s^{\Gw}(x,0)}=k. \ee
Take $\ge>0$ and put $u_\ge:=(1+\ge)u'+\ge$, $v:=(u-u_\ge)_+$. Then by \eqref{luba} there exists a smooth bounded domain $G \sbs \Gw$ such that $v=0$ in $G^c$ and $\tr_s^G(v)=0$. In light of Kato's inequality, we derive $\fw v \leq 0$ in $G$. Moreover, $v \leq kM_s^{\Gw}(\cdot,0)$ in $G$. By \rlemma{sub} we obtain $v =0$ in $G$ and therefore $u \leq (1+\ge)u' + \ge$ in $\Gw$. Letting $\ge \to 0$ yields $u \leq u'$ in $\Gw$. By permuting the role of $u$ and $u'$, we derive $u=u'$ in $\Gw$. 

By a similar argument as in step 2, we can show that $k \mapsto u_k^\Gw$ is increasing. \qeda \medskip

\noindent \textbf{Proof of Theorem F}.  (i) \textit{Case 1: $p_1^*<p<p_2^*$}.

Since $\prt \Gw \in C^2$, there exist two open balls $B$ and $B'$ such that $B\subset \Gw\subset B'^c$ and $\prt B \cap \prt B' = \{0 \}$. Since $M_s^B(x,0)\leq M_s^\Gw(x,0)\leq M_s^{B'^c}(x,0)$ it follows from \rth{unbduk} that \bel{uniq8} 
u_{k}^{B}\leq u_{k}^{\Gw}\leq u_{k}^{B'^c}
 \ee
where the first inequality holds in $B$ and the second inequality holds in $\Gw$. 

Let $\CO$ be $B$, $\Gw$ or $B'^c$. Because of uniqueness, we have
 \bel{uniq9} 
T_\ell[u^\CO_{k}]=u^{\CO_{\ell}}_{k\ell^{\frac{2s}{p-1}+1-N}} \qq\forall\ell>0,
 \ee
 with $\CO_{\ell}=\ell^{-1}\CO$. By \rth{unbduk}, the sequence $\{ u^\CO_{k}\}$ is increasing and by \eqref{ukU}, $u^\CO_{k} \leq U$. It follows that $\{u^\CO_{k}\}$ converges to a function  $u^\CO_{\infty}$ which is a positive solution of \eqref{eqa} with $\Gw$ replaced by $\CO$.  \smallskip
 
\noindent{\it Step 1: $\CO:=\BBR^N_+$}. Then $\CO_{\ell}=\BBR_+^N$. Letting $k \to\infty$ in \eqref{uniq9} yields to
 \bel{uniq10} 
T_\ell[u_{\infty}^{\BBR^N_+}]=u_{\infty}^{\BBR^N_+}\qq\forall\ell>0.
 \ee
Therefore $u_{\infty}^{\BBR^N_+}$ is self-similar and thus it can be written in the separable form 
$$ u_{\infty}^{\BBR^N_+}(x)=u_\infty^{\BBR_+^N}(r,\gs)=r^{-\frac{2s}{p-1}}\gw(\gs)
$$
where $r=|x|$, $\gs=\frac{x}{|x|} \in S^{N-1}$ and $\gw$ satisfies \eqref{I2}. Since $p_1^*<p<p_2^*$, it follows from Theorem E that $\gw=\gw^*$, the unique positive solution of \eqref{I2}.  This means
\bel{uniq11}
u_{\infty}^{\BBR^N_+}(x)=r^{-\frac{2s}{p-1}}\gw^*(\gs).
 \ee
This implies \eqref{uinf2}. \smallskip 
 
 \noindent{\it Step 2: $\CO:=B$ or $B'^c$}. In accordance with our previous notations, we set
 $B_{\ell}=\ell^{-1}B$ and $(B'^c)_\ell=\ell^{-1}B'^c$ for $\ell>0$ and we have,
 \bel{uniq11*} T_\ell[u_{\infty}^{B}]=u_{\infty}^{B_{\ell}} \text{ and }T_\ell[u_{\infty}^{B'^c}]=u_{\infty}^{(B'^c)_\ell}
 \ee
 and
  \bel{uniq12*} u_{\infty}^{B_{\ell'}}\leq u_{\infty}^{B_{\ell}} \leq u_{\infty}^{\BBR^N_+}\leq  u_{\infty}^{(B'^c)_\ell} 
 \leq  u_{\infty}^{(B'^c)_{\ell''}}\qq 0<\ell\leq\ell',\ell''\leq 1.
\ee
When $\ell\to 0$, $u_{\infty}^{B_{\ell}}\uparrow \underline u_{\infty}^{\BBR^N_+}$ and $u_{\infty}^{(B'^c)_\ell}\downarrow \overline u_{\infty}^{\BBR^N_+}$ where $ \underline u_{\infty}^{\BBR^N_+}$ and $\overline u_{\infty}^{\BBR^N_+}$ are positive solutions of \eqref{ukU} in $\BBR^N_+$ such that 
  \bel{uniq13} 
 u_{\infty}^{B_{\ell}} \leq \underline u_{\infty}^{\BBR^N_+}\leq u_{\infty}^{\BBR^N_+}\leq \overline u_{\infty}^{\BBR^N_+}\leq  u_{\infty}^{(B'^c)_\ell}\qq 0<\ell\leq 1.
\ee
Furthermore there also holds for $\ell,\ell'>0$,
  \bel{uniq14} 
T_{\ell'\ell}[u_{\infty}^{B}]=T_{\ell'}[T_{\ell}[u_{\infty}^{B}]]=u_{\infty}^{B_{\ell\ell'}} \text{ and }
T_{\ell'\ell}[u_{\infty}^{B'^c}]=T_{\ell'}[T_{\ell}[u_{\infty}^{B'^c}]]=u_{\infty}^{(B'_c)_{\ell\ell'}}.
\ee
Letting $\ell\to 0$ and using \eqref{uniq11*} and the above convergence, we obtain
  \bel{uniq15} 
\underline u_{\infty}^{\BBR^N_+}=T_{\ell'}[\underline u_{\infty}^{\BBR^N_+}] \text{ and }
\overline u_{\infty}^{\BBR^N_+}=T_{\ell'}[\overline u_{\infty}^{\BBR^N_+}] \qquad \ell'>0.
\ee
Again this implies that $\underline u_{\infty}^{\BBR^N_+}$ and $\overline u_{\infty}^{\BBR^N_+}$ are separable solutions of \eqref{I1}. Since $p_1^*<p<p_2^*$, by Theorem E, 
$$\underline u_{\infty}^{\BBR^N_+}(x)=\overline u_{\infty}^{\BBR^N_+}(x)=u_{\infty}^{\BBR^N_+}(x)=r^{-\frac{2s}{p-1}}\gw^*(\gs) \quad \text{with } r=|x|, \;\; \gs=\frac{x}{|x|}, \quad x \neq 0.$$
\smallskip 
 
 \noindent{\it Step 3: End of the proof}. From \eqref{uniq8} and \eqref{uniq11*} there holds
   \bel{uniq16} u_{\infty}^{B_{\ell}}\leq T_{\ell}[u_{\infty}^{\Gw}]
 \leq  u_{\infty}^{(B'^c)_\ell}\qq 0<\ell\leq 1.
\ee
Since the left-hand side and the right-hand side of \eqref{uniq16} converge to the same function $u_{\infty}^{\BBR^N_+}$, we obtain 
   \bel{uniq17} 
\lim_{\ell\to 0}\ell^{\frac{2s}{p-1}}u_{\infty}^{\Gw}(\ell x)=|x|^{-\frac{2s}{p-1}}\gw^*(\frac{x}{|x|})
\ee
and this convergence holds in any compact subset of $\Gw$. Take $|x|=1$, we derive \eqref{min1}. Estimate \eqref{uinf2} follows from \rprop{uinf1}.

(ii) \textit{Case 2: $0<p\leq p_1^*$.} Then by \rprop{uinf2}, $\lim_{k \to \infty}u_k^\Gw(x)=\infty$ for every $x \in \Gw$.  
\qeda

\appendix

\section{Appendix - Separable solutions}

\subsection{Separable $s$-harmonic functions}
We denote by $(r,\gs)\in \BBR_+\ti S^{N-1}$ the spherical  coordinates in $\BBR^N$, consider the following parametric representation of the unit sphere
\bel{S--1}\BA {lll}\displaystyle
S^{N-1}=\left\{\gs=(\cos\gf\,\gs',\sin\gf):\gs'\in S^{N-2}, -\tfrac\gp2\leq\gf\leq \tfrac\gp2\right\},
\EA\ee
hence $x_N=r\sin\gf$. We define the spherical fractional Laplace-Beltrami operator $\CA_{s}$ by
\bel{S-0}\BA {lll}\displaystyle
\CA_{s}\gw(\gs):=\lim_{\ge\to 0}\CA_{s,\ge}\gw(\gs)
\EA\ee
with
\bel{S-1}\BA {lll}\displaystyle
\CA_{s,\ge}\gw(\gs):=a_{N,s}\myint{}{}\!\!\myint{\BBR_+\ti S^{N-1}\setminus B_\ge(\overrightarrow\gs)}{}
\myfrac{(\gw(\gs)-\gw(\eta))\gt^{N-1}}{\left(1+\gt^2-2\gt\langle\gs,\eta\rangle\right)^{\frac{N}{2}+s}}dS(\eta)d\gt
\EA\ee
where $\overrightarrow\gs=(1,\gs)$. If $u:(r,\gs)\mapsto u(r,\gs)=r^{-\gb}\gw(\gs)$ is $s$-harmonic in $\BBR^N\setminus\{0\}$, it satisfies, at least formally, 
\bel{T-0}\BA {lll}\displaystyle
\CA_{s}\gw-\CL_{s,\gb}\gw=0\qquad\text{on }\,S^{N-1}
\EA\ee
where $\CL_{s,\gb}$ is the integral operator  
\bel{S-2}\BA {lll}\displaystyle
\CL_{s,\gb}\gw(\gs):=a_{N,s}\myint{0}{\infty}\myint{S^{N-1}}{}
\myfrac{(\gt^{-\gb}-1)\gt^{N-1}}{\left(1+\gt^2-2\gt\langle\gs,\eta\rangle\right)^{\frac{N}{2}+s}}\gw(\eta) dS(\eta)d\gt,
\EA\ee
whenever this integral is defined. We will see in the next two lemmas that the role of the exponent $\gb_0=N$ is fundamental for the definition of $\CL_{s,\gb}\gw$ since we have
\blemma{CL1} If $N\geq 2$, $s\in (0,1)$, $\gb<N$ and $(\gs,\eta)\in \BBR^{N-1}\ti \BBR^{N-1}$ such that $\langle\gs,\eta\rangle\neq 1$, we define
\bel{S-3}\BA {lll}\displaystyle
B_{s,\gb}(\gs,\eta):=\myint{0}{\infty}\myfrac{(\gt^{-\gb}-1)\gt^{N-1}}{\left(1+\gt^2-2\gt\langle\gs,\eta\rangle\right)^{\frac{N}{2}+s}}d\gt.
\EA\ee
Then \smallskip

\nind (i)\phantom{ii} $B_{s,\gb}(\gs,\eta)<0\Longleftrightarrow \gb<N-2s$,\smallskip

\nind (ii)\phantom{i} $B_{s,\gb}(\gs,\eta)=0\Longleftrightarrow \gb=N-2s$,\smallskip

\nind (iii) $B_{s,\gb}(\gs,\eta)>0\Longleftrightarrow \gb>N-2s$.
\es

\nind\proof Since $\gb<N$, the integral in \eqref{S-4} is absolutely convergent. We write
$$\BA {lll}B_{s,\gb}(\gs,\eta)=\myint{0}{1}\myfrac{(\gt^{-\gb}-1)\gt^{N-1}}{\left(1+\gt^2-2\gt\langle\gs,\eta\rangle\right)^{\frac{N}{2}+s}}d\gt+\myint{1}{\infty}\myfrac{(\gt^{-\gb}-1)\gt^{N-1}}{\left(1+\gt^2-2\gt\langle\gs,\eta\rangle\right)^{\frac{N}{2}+s}}d\gt\\[2mm]
\phantom{B_{s,\gb}(\gs,\eta)}=:I+II.
\EA$$
By the change of variable $\gt\mapsto\gt^{-1}$
$$II=-\myint{0}{1}\myfrac{(\gt^{-\gb}-1)\gt^{N-1+c_s}}{\left(1+\gt^2-2\gt\langle\gs,\eta\rangle\right)^{\frac{N}{2}+s}}d\gt,
$$
where
$c_s=\gb+2s-N$. Since
\bel{S-4}B_{s,\gb}(\gs,\eta)=\myint{0}{1}\myfrac{(\gt^{-\gb}-1)(\gt^{N-1}-\gt^{N-1+c_s})}{\left(1+\gt^2-2\gt\langle\gs,\eta\rangle\right)^{\frac{N}{2}+s}}d\gt,
\ee
the claim follows.\qeda\medskip

 As a byproduct of \eqref{S-4} we have the following monotonicity formula
 
 \blemma{CL1-1} If $N\geq 2$ and $s\in (0,1)$, then for any $(\gs,\eta)\in S^{N-1}\ti S^{N-1}$ the mapping $\gb\mapsto B_{s,\gb}(\gs,\eta)$ is continuous and increasing from $(N-2s,N)$ onto $(0,\infty)$.
\es

In the next result we analyze the behavior of $B_{s,\gb}(\gs,\eta)$ when $\gs-\eta\to 0$ on $S^{N-1}$.

\blemma{CL2} Assume $N\geq 2$, $s\in (0,1)$ and $\gb<N$ with $\gb\neq N-2s$, then \smallskip

\nind I- If $N\geq 3$, there exists $c=c(N,\gb,s)>0$ such that 
\bel{S-5}\BA {lll}\displaystyle
\abs{B_{s,\gb}(\gs,\eta)}\leq c\abs{\gs-\eta}^{3-N-2s} \qquad\forall (\gs,\eta)\in S^{N-1}\ti S^{N-1}.
\EA\ee
\nind II- If $N=2$, \\
(i) either $s>\frac{1}{2}$ and \eqref{S-5} holds with $N=2$,\\
(ii) either $s=\frac{1}{2}$ and 
\bel{S-6}\BA {lll}\displaystyle
\abs{B_{s,\gb}(\gs,\eta)}\leq c\left(-\ln\abs{\gs-\eta} +1\right)  \qquad\forall (\gs,\eta)\in S^{1}\ti S^{1}
\EA\ee
(iii) or $0<s<\frac{1}{2}$ and 
\bel{S-7}\BA {lll}\displaystyle
\abs{B_{s,\gb}(\gs,\eta)}\leq c  \qquad\forall (\gs,\eta)\in S^{1}\ti S^{1}
\EA\ee
\es

\nind\proof First, notice that the quantity 
$$\myint{0}{\frac{1}{2}}\myfrac{(\gt^{-\gb}-1)(\gt^{N-1}-\gt^{N-1+c_s})}{\left(1+\gt^2-2\gt\langle\gs,\eta\rangle\right)^{\frac{N}{2}+s}}d\gt
$$
is uniformly bounded with respect to $(\gs,\eta)$. The only possible singularity in the expression given in \eqref{S-4} occurs when $\langle\gs,\eta\rangle =1$ and $\gt=1$. 
We write $\langle\gs,\eta\rangle =1-\frac12\gk^2$ and $t=1-\gt$, hence
$$\BA {lll}\left(1+\gt^2-2\gt\langle\gs,\eta\rangle\right)^{\frac{N}{2}+s}
=\left(t^2+(1-t)\gk^2\right)^{\frac{N}{2}+s}\\[1mm]
\phantom{\left(1+\gt^2-2\gt\langle\gs,\eta\rangle\right)^{\frac{N}{2}+s}}
\approx \gk^{N+2s}\left(1+\left(\frac{t}{\gk}\right)^2\right)^{\frac{N}{2}+s}
\EA$$
as $t\to 0$. Moreover
$$\BA {lll}(\gt^{-\gb}-1)(\gt^{N-1}-\gt^{N-1+c_s})=
((1-t)^{-\gb}-1)((1-t)^{N-1}-(1-t)^{N-1+c_s})\\[2mm]
\phantom{(\gt^{-\gb}-1)(\gt^{N-1+c_s}-\gt^{N-1})}
=c_s\gb t^2+O(t^3)\quad\text{as }\;t\to 0.
\EA$$
Hence
$$\BA {lll}\myint{\frac{1}{2}}{1}\myfrac{(\gt^{-\gb}-1)(\gt^{N-1}-\gt^{N-1+c_s})}{\left(1+\gt^2-2\gt\langle\gs,\eta\rangle\right)^{\frac{N}{2}+s}}d\gt=
\myint{0}{\frac{1}{2}}\myfrac{((1-t)^{-\gb}-1)((1-t)^{N-1}-(1-t)^{N-1+c_s})}{\left(t^2+(1-t)\gk^2\right)^{\frac{N}{2}+s}}dt\\[4mm]
\phantom{\myint{\frac{1}{2}}{1}\myfrac{(\gt^{-\gb}-1)(\gt^{N-1+c_s}-\gt^{N-1})}{\left(1+\gt^2-2\gt\langle\gs,\eta\rangle\right)^{\frac{N}{2}+s}}d\gt}
\approx c_s\gk^{3-N-2s}\myint{0}{\frac{1}{2\gk}}\myfrac{x^2}{(1+x^2)^{\frac{N}{2}+s}}dx.
\EA$$
If $N=2$ and $s<\frac{1}{2}$, 
$$\abs{\gk^{1-2s}\myint{0}{\frac{1}{2\gk}}\myfrac{x^2}{(1+x^2)^{1+s}}dx}\leq M
$$
for some $M>0$ independent of $\gk$. If  $N=2$ and $s=\frac{1}{2}$
$$\myint{0}{\frac{1}{2\gk}}\myfrac{x^2}{(1+x^2)^{1+\frac{1}{2}}}dx= \ln \left(\frac{1}{\gk}\right) (1+o(1))
$$
and if $N=3$ or $N=2$ and $s>\frac{1}{2}$,
$$\myint{0}{\frac{1}{2\gk}}\myfrac{x^2}{(1+x^2)^{\frac{N}{2}+s}}dx\to \myint{0}{\infty}\myfrac{x^2}{(1+x^2)^{\frac{N}{2}+s}}dx.
$$
as $\gk\to 0$. Since $\gs,\eta \in S^{N-1}$ there holds $\gk^2=2(1-\langle\gs,\eta\rangle)=\abs{\gs-\eta}^2$. Thus the claim follows.\qeda


\bprop{CL} Assume $N\geq 2$, $s\in (0,1)$ and $\gb<N$ with $\gb\neq N-2s$. Then $\gw\mapsto \CL_{s,\gb}\gw$ is a continuous linear operator from $L^q(S^{N-1})$ into $L^r(S^{N-1})$ for any $1\leq q, r\leq\infty$ such that 
\bel{S-8}\BA {lll}\displaystyle
\myfrac{1}{r}>\myfrac{1}{q}-\myfrac{2(1-s)}{N-1}.
\EA\ee
Furthermore, $ \CL_{s,\gb}$ is positive (resp. negative) operator  if $\gb<N-2s$ (resp. $N-2s<\gb<N$).
\es
\nind\proof By \rlemma{CL2}, for any $\eta \in S^{N-1}$,  $B_{s, \gb}(.,\eta)\in L^a(S^{N-1})$ for all $1<a<\frac{N-1}{N+2s-3}$ if $N\geq 3$ or $N=2$ and $s>\frac{1}{2}$; $ B_{s,\gb}(.,\eta)\in\bigcap_{1\leq a<\infty}L^a(S^{1})$ if $N=2$ and $s=\frac{1}{2}$ and $ B_{s,\gb}(.,\eta)$ is uniformly bounded on $S^{1}$ if $N=2$ and $0<s<\frac{1}{2}$. The continuity result follows from Young's inequality and the sign assertion from 
\rlemma{CL1}. \qeda\medskip

The above calculations justifies the name of fractional Laplace-Beltrami operator given to $\CA_{s}$ since we have  the following relation.

\blemma{CL0} Assume $N\geq 2$ and $s\in (0,1)$, then
\bel{S-1-1}\BA {lll}\displaystyle
\CA_{s}\gw(\gs)=b_{N,s}CPV\myint{S^{N-1}}{}\myfrac{(\gw(\gs)-\gw(\eta))}{\abs{\gs-\eta}^{N-1+2s}}dS(\eta) +\CB_{s}\gw(\gs), 
\EA\ee
where $\CB_{s}$ is a bounded linear operator from $L^q(S^{N-1})$ into $L^r(S^{N-1})$ for $q$, $r$ satisfying \eqref{S-8} and 
\bel{S-1-2}\BA {lll}\displaystyle
b_{N,s}:=2a_{N,s}\myint{0}{\infty}\myfrac{dx}{(x^2+1)^{\frac{N}{2}+s}}.
\EA\ee
\es
\nind\proof If $(\gs,\eta)\in S^{N-1}\ti S^{N-1}$, we set $\langle\gs,\eta\rangle=1-\frac12\gk^2$. Then
$$\myint{0}{\infty}\myfrac{\gt^{N-1}d\gt}{\left(1+\gt^2-2\gt\langle\gs,\eta\rangle\right)^{\frac{N}{2}+s}}
=\myint{0}{1}\myfrac{\left(\gt^{N-1}+\gt^{2s-1}\right)d\gt}{\left(1+\gt^2-2\gt\langle\gs,\eta\rangle\right)^{\frac{N}{2}+s}}.
$$
Then we put $t=1-\gt$, hence, when $t\to 0$, we have after some straightforward computation
$$\BA {lll}\myfrac{\left(\gt^{N-1}+\gt^{2s-1}\right)}{\left(1+\gt^2-2\gt\langle\gs,\eta\rangle\right)^{\frac{N}{2}+s}}
=\myfrac{\left(2-(N+2s -2)t+O(t^2)\right)\left(1+\frac{(N+2s)t\gk^2}{  2(t^2+\gk^2)}+O\left(\left(\frac{t\gk^2}{t^2+2\gk^2}\right)^2\right)\right)}{(t^2+\gk^2)^{\frac{N}{2}+s}}\\[4mm]
\phantom{\myfrac{\left(\gt^{N-1}+\gt^{2s-1}\right)}{\left(1+\gt^2-2\gt\langle\gs,\eta\rangle\right)^{\frac{N}{2}+s}}}
=\myfrac{2+ 2t + O(t^2)}{(t^2+\gk^2)^{\frac{N}{2}+s}}.
\EA$$
This implies
\begin{align} \BAL &\myint{0}{1}\myfrac{\left(\gt^{N-1}+\gt^{2s-1}\right)d\gt}{\left(1+\gt^2-2\gt\langle\gs,\eta\rangle\right)^{\frac{N}{2}+s}}\\
&=2\gk^{1-N-2s}\myint{0}{\frac 1\gk}\myfrac{dx}{(x^2+1)^{\frac{N}{2}+s}}
+  2\gk^{2-N-2s}\myint{0}{\frac 1\gk}\myfrac{x dx}{(x^2+1)^{\frac{N}{2}+s}} 
+O(\gk^{3-N-s})\myint{0}{\frac 1\gk}\myfrac{ x^2  dx}{(x^2+1)^{\frac{N}{2}+s}}\\
&=2\gk^{1-N-2s}\myint{0}{\infty}\myfrac{dx}{(x^2+1)^{\frac{N}{2}+s}}+O(1)+O(\gk^{3-N-s})\myint{0}{\frac 1\gk}\myfrac{  x^2 dx}{(x^2+1)^{\frac{N}{2}+s}}. 
\EAL \end{align}
Since $\gk=\abs{\gs-\eta}$, the claim follows from \rprop{CL} and the kernel estimate in \rlemma{CL2}.\qeda

\blemma{CL3} Under the assumption of \rlemma{CL0} there holds
\bel{S-9}\BA {lll}\displaystyle
\abs{\myint{S^{N-1}}{}\gw\CL_{s,\gb}\gw dS}\leq c_{35}\myint{S^{N-1}}{}\gw^2dS\qquad\forall\gw\in L^{2}(S^{N-1}),
\EA\ee
where
$$c_{35}=\myint{0}{1}\left(\myint{S^{N-1}}{}\myfrac{dS(\eta)}{(1+\gt^2-2\gt\langle{{\bf e}_{_N}},\eta\rangle)^{\frac{N}{2}+s}}\right)(\gt^{-\gb}-1)\abs{\gt^{N-1}-\gt^{N-1+c_s}}d\gt.
$$
\es
\nind\proof There holds by Cauchy-Schwarz inequality
$$\BAL &\abs{\myint{S^{N-1}}{}\gw\CL_{s,\gb}\gw dS}\\
&\leq \myint{0}{1}\!\left(\myint{S^{N-1}}{}\!\myint{S^{N-1}}{}\!\myfrac{|\gw(\eta)||\gw(\gs)|dS(\eta)dS(\gs)}{\left(1+\gt^2-2\gt\langle\gs,\eta\rangle\right)^{\frac{N}{2}+s}}\right)\!(\gt^{-\gb}-1)\abs{\gt^{N-1}-\gt^{N-1+c_s}}d\gt\\
&\leq \myint{0}{1}\!\left(\myint{S^{N-1}}{}\myint{S^{N-1}}{}\myfrac{\gw^2(\eta)}{\left(1+\gt^2-2\gt\langle\gs,\eta\rangle\right)^{\frac{N}{2}+s}}dS(\eta)dS(\gs)\right)\ti (\gt^{-\gb}-1)\abs{\gt^{N-1}-\gt^{N-1+c_s}}d\gt
\\
&\leq \myint{S^{N-1}}{}\left(\myint{0}{1}\left(\myint{S^{N-1}}{}\myfrac{dS(\gs)}{\left(1+\gt^2-2\gt\langle\gs,\eta\rangle\right)^{\frac{N}{2}+s}}\right)(\gt^{-\gb}-1)\abs{\gt^{N-1}-\gt^{N-1+c_s}}d\gt\right)
\gw^2(\eta)dS(\eta).
\EAL $$
Since, by invariance by rotation, we have
$$\myint{S^{N-1}}{}\myfrac{dS(\gs)}{\left(1+\gt^2-2\gt\langle\gs,\eta\rangle\right)^{\frac{N}{2}+s}}
=\myint{S^{N-1}}{}\myfrac{dS(\gs)}{\left(1+\gt^2-2\gt\langle{{\bf e}_{_N}},\gs\rangle\right)^{\frac{N}{2}+s}},
$$
we derive \eqref{S-9}.\qeda \medskip

We denote the upper hemisphere of the unit sphere in $\BBR^N$ by $S_+^{N-1}=S^{N-1} \cap \BBR_+^N$.
\bprop{eig}Let $N\geq 2$, $s\in (0,1)$ and $N-2s<\gb<N$. Then there exist a unique $\gl_{s,\gb}>0$ and 
a unique (up to an homothety) positive $\psi_1\in W^{s,2}_0(S^{N-1}_+)$, such that 
\bel{S-13}\BA {lll}
\CA_s\psi_1=\gl_{s,\gb}\CL_{s,\gb}\psi_1\qquad\text {in }\,S^{N-1}_+.
\EA\ee
Furthermore the mapping $\gb\mapsto \gl_{s,\gb}$ is continuous and decreasing from $(N-2s,N)$ onto 
$(0,\infty)$. Finally $\gl_{s,\gb}=1$ if and only if $\gb=N-s$ and $\psi_1(\gs)=(\sin\gf)^s$.
\es

\nind\proof We first notice that 
\bel{S-14}\BA {lll}\displaystyle
\myint{S^{N-1}_+}{}\gw\CA_s\gw dS=\myfrac{1}{2}\myint{S^{N-1}_+}{}\myint{0}{\infty}\myint{S^{N-1}_+}{}\myfrac{\left(\gw(\gs)-\gw(\eta)\right)^2}{\left(1+\gt^2-2\gt\langle\gs,\eta\rangle\right)^{\frac{N}{2}+s}}\gt^{N-1}dS(\eta)d\gt d S(\gs),
 \EA\ee
for any $\gw\in C^1_0(S^{N-1}_+)$. By \rlemma{CL0} and \eqref{S-8} with $r=q=2$,
$$\BA {lll}
\myint{S^{N-1}_+}{}\myint{0}{\infty}\myint{S^{N-1}_+}{}\myfrac{\left(\gw(\gs)-\gw(\eta)\right)^2}{\left(1+\gt^2-2\gt\langle\gs,\eta\rangle\right)^{\frac{N}{2}+s}}\gt^{N-1}dS(\eta)d\gt d S(\gs)\\[4mm]
\phantom{---------------------}
\leq c_{36}\norm{ \gw}^2_{W^{s,2}_0(S^{N-1}_+)}+c_{37}\norm\gw^2_{L^2(S^{N-1}_+)},
 \EA$$
 where 
 $$\norm{ \gw}^2_{W^{s,2}_0(S^{N-1}_+)}=\myint{S^{N-1}_+}{}\myint{S^{N-1}_+}{}\myfrac{\left(\gw(\gs)-\gw(\eta)\right)^2}{\abs{\eta-\gs}^{N-1+2s}}dS(\eta) d S(\gs).
 $$
 Since, by Poincar\'e inequality \cite{NPV}, there holds
$$\norm{ \gw}^2_{W^{s,2}_0(S^{N-1}_+)}\geq c_{38}\norm{ \gw}^2_{L^2(S^{N-1}_+)},
$$
we obtain that the right-hand side of \eqref{S-14} is bounded from above by $\left(\frac12 c_{36}+\frac{c_{37}}{2c_{38}}\right)\norm{ \gw}^2_{W^{s,2}_0(S^{N-1}_+)}$. Next we use the expansion estimates in \rlemma{CL0} to obtain that 
$$\myfrac{\gt^{N-1}+\gt^{2s-1}}{\left(1+\gt^2-2\gt\langle\gs,\eta\rangle\right)^{\frac{N}{2}+s}}
\geq \myfrac{1}{(t^2+\gk^2)^{\frac{N}{2}+s}}\qquad \forall t=1-\gt\in(0, \ge_0)\,,\;\forall (\gs,\eta)\in S^{N-1}_+\ti S^{N-1}_+,
$$
where $\gk=\abs{\gs-\eta}\leq 2$. Hence
$$\BA {lll}\myint{0}{\infty}\myfrac{\gt^{N-1}d\gt}{\left(1+\gt^2-2\gt\langle\gs,\eta\rangle\right)^{\frac{N}{2}+s}}
\geq \myint{0}{\ge_0}\myfrac{dt}{(t^2+\gk^2)^{\frac{N}{2}+s}}=\gk^{1-N-2s}\myint{0}{\frac{\ge_0}{2}}\myfrac{dt}{(t^2+1)^{\frac{N}{2}+s}}.
\EA$$
Therefore, 
$$\myint{S^{N-1}_+}{}\gw\CA_s\gw dS\geq \myint{0}{\frac{\ge_0}{2}}\myfrac{dt}{2(t^2+1)^{\frac{N}{2}+s}}\norm{ \gw}^2_{W^{s,2}_0(S^{N-1}_+)}.
$$
Finally we obtain 
\begin{equation} \label{A.18} \BA {lll}
\myfrac{1}{c_{39}}\norm{ \gw}^2_{W^{s,2}_0(S^{N-1}_+)}\\[4mm]\phantom{--}
\leq\myint{S^{N-1}_+}{}\myint{0}{\infty}\myint{S^{N-1}_+}{}\myfrac{\left(\gw(\gs)-\gw(\eta)\right)^2}{\left(1+\gt^2-2\gt\langle\gs,\eta\rangle\right)^{\frac{N}{2}+s}}\gt^{N-1}dS(\eta)d\gt d S(\gs)\\[4mm]
\phantom{--------------------------}
\leq c_{39}\norm{ \gw}^2_{W^{s,2}_0(S^{N-1}_+)}.
 \EA \end{equation}
We consider the bilinear form in $W^{s,2}_0(S^{N-1}_+)$
$$\BBA(\gw,\gz):=\myint{S^{N-1}_+}{}\myint{0}{\infty}\myint{S^{N-1}_+}{}\myfrac{\left(\gw(\gs)-\gw(\eta)\right)\gz(\gs)}{\left(1+\gt^2-2\gt\langle\gs,\eta\rangle\right)^{\frac{N}{2}+s}}\gt^{N-1}dS(\eta)d\gt d S(\gs).
$$
Then $\BBA$ is symmetric and there holds
$$\BBA(\gw,\gw)=\myint{S^{N-1}_+}{}\gw\CA_s\gw dS\geq \myfrac{1}{2c_{39}}\norm{ \gw}^2_{W^{s,2}_0(S^{N-1}_+)},
$$
and 
$$\abs{\BBA(\gw,\gz)}\leq \left(\myint{S^{N-1}_+}{}\gw\CA_s\gw dS\right)^{\frac{1}{2}}\left(\myint{S^{N-1}_+}{}\gz\CA_s\gz dS\right)^{\frac{1}{2}}\leq \myfrac{c_{39}}{2}\norm{ \gw}_{W^{s,2}_0(S^{N-1}_+)}\norm{ \gz}_{W^{s,2}_0(S^{N-1}_+)}.
$$
By Riesz theorem, for any $L\in W^{-s,2}(S^{N-1}_+)$ there exists $\gw_L\in W^{s,2}_0(S^{N-1}_+)$ such that 
$$\BBA(\gw_L,\gz)=L(\gz)\qquad\forall\gz\in W^{s,2}_0(S^{N-1}_+).$$
We denote $\gw_L=\CA_s^{-1}(L)$. It is clear that $\CA_s^{-1}$ is positive and since the
 the embedding of $W^{s,2}_0(S^{N-1}_+)$ into $L^2(S^{N-1}_+)$ is compact by Rellich-Kondrachov theorem \cite{NPV}, $\CA_s^{-1}$ is a compact operator. Hence the operator
 $$\gw\mapsto \CA_s^{-1}\circ \CL_{s,\gb}\gw
 $$
is a compact positive operator (here we use the fact that $\gb>N-2s$ which makes $\CB_{s,\gb}$ positive). By the Krein-Rutman theorem there exists $\gm>0$ and $\psi_1\in W^{s,2}_0(S^{N-1}_+)$, 
$\psi_1\geq 0$ such that 
 $$\CA_s^{-1}\circ \CL_{s,\gb}\psi_1=\gm\psi_1.
 $$
The function $\psi_1$ is the unique positive eigenfunction and $\gm$ is the only positive eigenvalue with positive eigenfunctions. Furthermore $\gm$ is the spectral radius of $\CA_s^{-1}\circ \CB_{s,\gb}$. If we set $\gl_{s,\gb}=\gm^{-1}$, we obtain \eqref{S-13}. It is also classical that $\gl_{s,\gb}$ can be defined by
\bel{S-16}\BA {lll}
\gl_{s,\gb}:=\inf\left\{\int_{S_+^{N-1}} \gw \CA_s \gw dS :\gw\in W^{s,2}_0(S^{N-1}_+),\gw\geq 0,
\myint{S^{N-1}_+}{}\gw\CL_{s,\gb}\gw dS=1\right\}.
 \EA\ee
Using \eqref{S-4}, \rlemma{CL1-1} and monotone convergence theorem, we derive that the mapping
$$\gb\mapsto \myint{S^{N-1}_+}{}\gw\CL_{s,\gb}\gw dS
$$ 
is increasing and continuous.  This implies that $\gb\mapsto\gl_{s,\gb}$ is decreasing and continuous. Since 
$\myint{S^{N-1}_+}{}\gw\CL_{s,\gb}\gw dS\to\infty$ when $\gb \uparrow N$, the expression \eqref{S-16} implies that 
$\gl_{s,\gb}\to 0$ when $\gb \uparrow N$. Next, if $\gw\geq 0$ is an element of  $W^{s,2}_0(S^{N-1}_+)$ such that $\myint{S^{N-1}_+}{}\gw\CL_{s,\gb}\gw dS=1$, we derive from Poincar\'e inequality \cite {NPV} and \eqref{S-9},
$$\norm{ \gw}^2_{W^{s,2}_0(S^{N-1}_+)}\geq c_{38}\norm{ \gw}^2_{L^2(S^{N-1}_+)}\geq \myfrac{c_{38}}{c_{35}}.
$$
Since $c_{35}\to 0$ when $\gb \downarrow N-2s$, we infer that $\displaystyle\lim_{\gb \to N-2s}\gl_{s,\gb}= \infty$. Consequently the mapping $\gb\mapsto\gl_{s,\gb}$ is a decreasing homeomorphism from $(N-2s, N)$ onto $(0,\infty)$ and there exists a unique $\gb_s\in (N-2s, N)$ such that $\gl_{s,\gb_s}=1$. 
The following expression of the Martin kernel in $\BBR^N_+$ is classical, 
\bel{MRN} M_s^{\BBR_+^N}(x,y)=c_{N,s}\,x_N^s|x-y|^{-N} \quad \forall x \in \BBR_+^N, y \in \prt \BBR_+^N, \ee
hence, if $y=0$, it is a separable singular $s$-harmonic function expressed  in spherical coordinates with $x=(r,\gs)$ by
$$M_s^{\BBR^N_+}((r,\gs),0)=c_{N,s}r^{s-N}(\sin\gf)^s.
$$
This means that the function $\gs\mapsto \gw(\gs)=(\sin\gf)^s$, which vanishes on $\overline{S^{N-1}_-}$ and belongs to  $W^{s,2}_0(S^{N-1}_+)\cap L^{\infty}(S^{N-1}_+)$, satisfies
$$\CA_s\gw-\CL_{s,N-s}\gw=0.
$$
The uniqueness of the positive eigenfunction implies that this function is $\psi_1$ and $\gb=N-s$.
\qeda
\medskip


\subsection{The nonlinear problem}
\subsubsection{Separable solutions in $\BBR^N$}
If we look for separable positive solutions of
\bel{S-17}\BA {lll}
(-\Delta)^su+u^p=0\qquad&\text{in }\;\BBR^N,
\EA\ee
under the form $u(x)=r^{-\frac{2s}{p-1}}\gw(\gs)$ where $x=(r,\gs)\in \BBR_+\ti S^{N-1}$, then
$\gw$ satisfies
\bel{S-18}\BA {lll}
\CA_s\gw-\CL_{s,\frac{2s}{p-1}}\gw+\gw^p=0\qquad&\text{in }\;S^{N-1}.
\EA\ee
\bprop{consta} Assume $N\geq 2$ and $s\in (0,1)$.  

(i) If $p \geq p_3^*$ then there exists no positive solution of \eqref{S-18}.

(ii) If $p_1^*<p<p_3^*$ then the unique positive solution of \eqref{S-18} is a constant function with value
\bel{S-19}\ell_{s,p}=\left(c_{35}\right)^{\frac{1}{p-1}},
\ee
where $c_{35}$ is the constant defined in \rlemma{CL3}.
\es 
\nind\proof   If $p\geq p_3^*$, we assume that there exists a solution $\gw\geq 0$ of \eqref{S-18}. Then $\gw$ satisfies
$$\myint{S^{N-1}}{}\gw\CA_s\gw dS-\myint{S^{N-1}}{}\gw\CL_{s,\frac{2s}{p-1}}\gw dS+\myint{S^{N-1}}{}\gw^{p+1}dS=0.
$$
Since $p \geq p_3^*$, we have $c_s\leq 0$ which implies 
$$\myint{S^{N-1}}{}\gw\CL_{s,\frac{2s}{p-1}}\gw dS\leq 0.$$ 
Then $\gw=0$ since the two other integrals are nonnegative. 

Next, if $p_1^*<p<p_3^*$ it is clear that if  $\gw$ is a constant nonnegative solution of \eqref{S-18} then we have  
$$\gw\myint{0}{1}\myint{S^{N-1}}{}\myfrac{(\gt^{-\frac{2s}{p-1}-1})(\gt^{N-1}-\gt^{N-1+c_s})}{(1+\gt^2-2\gt\langle\gs,\eta\rangle)^{\frac{N}{2}+s}}dS(\eta)d\gt=\gw^p\qquad\forall\gs\in S^{N-1}.
$$
Using invariance by rotation of the integral term on $S^{N-1}$, we derive the claim. Conversely, assume $\gw$ is any bounded nonconstant positive solution, then it belongs to $C^2(S^{N-1})$ by \cite{RS}. Let $\gs_0\in S^{N-1}$ where 
$\gw$ is maximal, then $\CA_s\gw(\gs_0)\geq 0$ thus
$$\gw^p(\gs_0)\leq \CL_{s,\frac{2s}{p-1}}\gw(\gs_0)\leq \gw(\gs_0)\myint{0}{1}\myint{S^{N-1}}{}
\myfrac{(\gt^{-\frac{2s}{p-1}}-1)(\gt^{N-1}-\gt^{N-1+c_s})}{(1+\gt^2-2\langle\gs_0,\eta)^{\frac{N}{2}+s}}dS(\eta) d\gt
=c_{35}\gw(\gs_0).$$
Hence $\gw(\gs_0)<\ell_{s,p}$. Similarly $\displaystyle\min_{S^{N-1}}\gw>\ell_{s,p}$, which is a contradiction. \qeda
\bcor{rad} Assume $N\geq 2$, $s\in (0,1)$ and $p_1^*<p<p_3^*$. Then the only positive separable solution 
$u$ of \eqref{S-17} in $\BBR^N\setminus\{0\}$ is 
\bel{S-20}
x\mapsto U(x)=\ell_{s,p}\abs x^{-\frac{2s}{p-1}}.
\ee
\es
\subsubsection{Separable solutions in $\BBR_+^N$}
If we consider separable solutions $x\mapsto u(x)=r^{-\frac{2s}{p-1}}\gw(\gs)$ of problem \eqref{I1} then $\gw$ satisfies \eqref{I2}.\smallskip

\noindent \textbf{Proof of Theorem E.}
 
\nind{\it Step 1: Non-existence}. Assume that such a solution $\gw\geq 0$ exists, then
$$\myint{S^{N-1}_+}{}\gw\CA_{s}\gw dS-\myint{S^{N-1}_+}{}\gw\CL_{s,\frac{2s}{p-1}}\gw dS+\myint{S^{N-1}_+}{}\gw^pdS=0.
$$
Hence
\bel{S1}
\left(\gl_{s,\frac{2s}{p-1}}-1\right)\myint{S^{N-1}_+}{}\gw\CL_{s,\frac{2s}{p-1}}\gw dS+\myint{S^{N-1}_+}{}\gw^pdS\leq 0.
\ee
If $\gl_{s,\frac{2s}{p-1}}\geq 1$, equivalently $p\geq p_2^*$, the only nonnegative solution is the trivial one.\smallskip

\nind{\it Step 2: Existence}. Consider the following functional with domain $W^{s,2}_0(S^{N-1}_+)\cap L^{p+1}(S^{N-1}_+)$,

\bel{S2} 
\gw\mapsto\CJ(\gw):=\int_{S_+^{N-1}}\gw \CA_s \gw dS +\myfrac{1}{p+1}\myint{S^{N-1}_+}{}\abs\gw^{p+1}dS-\myint{S^{N-1}_+}{}\gw\CL_{s,\frac{2s}{p-1}}\gw dS.
\ee
Because of \rlemma{CL3}, $\CJ(\gw)\to\infty$ when 
$\norm{ \gw}_{W^{s,2}_0(S^{N-1}_+)}+\norm{ \gw}_{L^{p+1}(S^{N-1}_+)}\to\infty$. Furthermore, for $\ge>0$, we have
$$\CJ(\ge\psi_1)=\ge^2\left(\gl_{s,\frac{2s}{p-1}}-1\right)\myint{S^{N-1}_+}{}\psi_1\CL_{s,\frac{2s}{p-1}}\psi_1 dS+\myfrac{\ge^{p+1}}{p+1}
\myint{S^{N-1}_+}{}\abs{\psi_1}^{p+1}dS.
$$
This implies that $\inf\CJ(\gw)<0$ if  $\gl_{s,\frac{2s}{p-1}}< 1$, and thus the infimum of $\CJ$ in $W^{s,2}_0(S^{N-1}_+)\cap L_+^{p+1}(S^{N-1}_+)$ is achieved by a nontrivial nonnegative solution of \eqref{I2}. 
\smallskip

\nind{\it Step 3: Uniqueness}. \smallskip

\nind{\it (i) Existence of a maximal solution.} By \cite{RS} any solution $\gw$ is smooth. Hence, 
at its maximum $\gs_0$, it satisfies $\CA_{s}\gw(\gs_0)\geq 0$, thus
$$\gw(\gs_0)^p\leq \CL_{s,\frac{2s}{p-1}}\gw(\gs_0)\leq \gw(\gs_0)c_{35}.
$$
This implies that $\sup\gw\leq \ell_{s,p}$. From the equation the set $\CE\subset W^{s,2}_0(S^{N-1}_+)$ of positive solutions of 
\eqref{I2} is bounded in $W^{s,2}_0(S^{N-1}_+)\cap L^{\infty}(S^{N-1}_+)$ and thus in $C^s(S^{N-1})\cap C^2(S^{N-1}_+)$ by \cite{RS}. We put $\overline\gw(\gs)=\sup\{\gw(\gs):\gw\in\CE\}$. There exists a countable dense set $\CS:=\{\gs_n\}\subset S^{N-1}_+$ and a sequence of function $\{\gw_n\}\subset \CE$ such that 
$$\lim_{n\to\infty}\gw_n(\gs_k)=\overline\gw(\gs_k).
$$
 Furthermore, this sequence $\{\gw_n\}$ can be constructed such that $\{\gw_n(\gs_k)\}$ is nondecreasing for any $k$. Finally by local compactness estimate, $\{\gw_n\}$ converges to $\overline\gw$ in $C^{s-\gd}(S^{N-1})\cap C^2(S^{N-1}_+)$
for any $\gd \in (0,s)$ and weakly in $W^{s,2}_0(S^{N-1}_+)$. This implies that $\overline\gw$ belongs to $\CE$. It follows from \cite[Th 1.2]{RS} that any $\gw\in \CE$ satisfies
\bel{S-23}\BA {lll}
\gw(\gs)\leq c_{40}\left(\dist(\gs,\prt S^{N-1}_+)\right)^{s}=c_{40}\gf^s\qquad\forall \gs \in S^{N-1}_+.
\EA\ee\smallskip

\nind{\it (ii) Existence of a minimal solution.} This follows from \rth{unbduk} that 
$u_k^{\BBR_+^N} \uparrow u_\infty^{\BBR_+^N}$ and $u_\infty^{\BBR_+^N}$ is self-similar and it is the minimal solution of \eqref{eqa} in $\BBR^N_+$ which satisfies 
\bel{S-23'}\BA {lll}\displaystyle
\lim_{x\to 0}\myfrac{u_\infty^{\BBR_+^N}(x)}{M_s^{\BBR^N_+}(x,0)}=\infty.
\EA\ee
Thus $u_\infty^{\BBR_+^N}(r,\gs)=r^{-\frac{2s}{p-1}}\underline\gw(\gs)$ and $\underline\gw$ is the minimal positive solution of \eqref{I2}.
Furthermore  it follows from \eqref{uinf2} that 
\bel{S-24}
\underline\gw(\gs)\geq c_{41}\left(\dist(\gs,\prt S^{N-1}_+)\right)^{s}=c_{41}\gf^s\qquad\forall \gs \in S^{N-1}_+,
\ee
if $\gf=\gf(\gs)$ is the latitude of $\gs$. \smallskip

\nind \nind{\it (iii) End of the uniqueness proof.} By combining \eqref{S-23} and \eqref{S-24} we infer that there exists $K>1$ such that 
\bel{S-25}\BA {lll}\displaystyle
\overline\gw\leq K\underline\gw\qquad\text{in }\,S^{N-1}_+.
\EA\ee
Assume $\overline\gw\neq\underline\gw$, then 
$$\gw_1:=\underline\gw-\myfrac{1}{2K}\left(\overline\gw-\underline\gw\right)
$$
is a positive supersolution (by convexity) of \eqref{I2}. Moreover 
$$\gw_2:=\left(\myfrac{1}{2}+\myfrac{1}{2K}\right)\underline\gw
$$
is a positive subsolution of \eqref{I2} smaller than $\gw_1$ hence also than $\underline\gw$. It follows by classical construction that there exists a solution $\tilde \gw$ of \eqref{I2} which satisfies $\gw_2\leq\tilde\gw\leq\gw_1$, which contradicts the minimality of $\underline\gw$. \qeda 

\bigskip

\noindent \textbf{Acknowledgements.} The first author is supported by Fondecyt Grant 3160207. The second author is supported by collaboration programs ECOS C14E08. 

\end{document}